\documentclass[10pt]{amsart}

\usepackage[edge-length=.5cm,root-radius=.05cm]{dynkin-diagrams}

\usepackage{hyperref}
\usepackage{breakurl}
\usepackage{graphics}
\usepackage{graphicx}

 \usepackage[normalem]{ulem}
\usepackage{amsfonts,amssymb,amsmath}

\usepackage{enumerate}
\usepackage{mathrsfs}
\usepackage{array,multirow, makecell}
\usepackage[all]{xy}

\usepackage[lofdepth,lotdepth]{subfig}
\newtheorem{theorem}{Theorem}
\newtheorem{proposition}[theorem]{Proposition}
\newtheorem{lemma}[theorem]{Lemma}% 
\theoremstyle{definition}
\newtheorem{example}{Example}%
\newtheorem{remark}{Remark}%
\newtheorem{notation}{Notation}

\newtheorem{definition}{Definition}%

\def\ZZ{\mathbb{Z}}
\def\PP{\mathbb{P}}
\def\NN{\mathbb{N}}
\def\QQ{\mathbb{Q}}

\newcommand{\Prim}[1]{\operatorname{Prim}(#1)}
\newcommand{\Dual}[2]{\langle #1,#2\rangle}
\newcommand{\Scalar}[2]{\left( #1,#2\right)}

\newcommand{\Conv}{\operatorname{Conv}}

\newcommand{\aff}{{\operatorname{aff}}}
\begin{document}

\title[Gorenstein Fano generic torus orbit closures in $G/P$]{Gorenstein Fano generic torus orbit closures in $\mathbf{G/P}$}


\author{Pierre-Louis Montagard}
\email{pierre-louis.montagard@umontpellier.fr}
\address{IMAG, Universit\'e de Montpellier, Place Eug\`ene Bataillon, Montpellier, 34090, France}

\author{Alvaro Rittatore}
\email{alvaro.rittatore@cmat.edu.uy}
\address{Centro de matem\'atica, Universidad de Montevideo, Igu\'a 4225 esq. Mataojo., Montevideo,Uruguay}

\begin{abstract}
	Given a reductive group \text{$\textstyle G$} and a parabolic
        subgroup \text{$\textstyle P\subset G$}, with maximal
	torus \text{$\textstyle T$}, we consider (following Dabrowski's work) the closure \text{$\textstyle X$} of a generic \text{$\textstyle T$}-orbit in \text{$\textstyle G/P$}, and determine in combinatorial terms
	when  the toric variety  \text{$\textstyle X$} is  $\QQ$-Gorenstein Fano, 
	extending  in this way the classification of smooth Fano
        generic closures given by Voskresenski\u{\i} and Klyachko. As an application,  we apply the well known
        correspondence between Gorenstein Fano toric varieties and
        reflexive polytopes in  order to exhibit which reflexive polytopes
        correspond to generic closures --- this list includes  the
        reflexive root polytopes.

\end{abstract}

	\maketitle
\section{Introduction}

 Toric varieties
 --- that is, normal varieties $X$ over an algebraically
 closed field $\Bbbk$, on which an algebraic 
 torus $T$ acts effectively and with an open orbit ---
 have been
 thoroughly studied since the beginning of the 1970's. Since the
 geometric properties of a toric variety can be described in
 combinatorial terms (by means of its associated fan), this family of
 algebraic varieties provides a nice framework in which to study
 either their geometric properties or the combinatorial properties of
 their associated fans.

In the early beginning of the theory of toric varieties, D.  Mumford
considered the toric variety associated to the fan 
obtained by considering the weight lattice  and the set of all closed
Weyl chambers of a root system $R$ (see 
\cite{kn:Kempf}). Afterwards, the geometry of this variety was  intensively 
studied by several authors (see \cite{kn:pro90,kn:DolgaLunts,kn:BriJoshua,kn:BatBlume}). In \cite{kn:voskly}, V.E.~Voskresenski\u{\i} and A.A.~Klyachko
considered a larger family of fans  constructed by fixing a set $I$ of
simple roots of $R$ and 
``gluing together'' selected adjacent Weyl chambers that
correspond to a choice of a proper subset  $L\subsetneq I$ (see
Definition \ref{def:FanRoot}).   
The invariance properties of these fans (w.r.t. the action of the Weyl group of $R$)
 allow the authors  to
characterize the pairs  $(R,L)$ such that the associated toric variety
$X_{R,L}$ is smooth Fano.    
A remarkable result of R.~Dabrowski proves that   the
toric varieties $X_{R,L}$ can be constructed as the closure of an 
orbit of a maximal torus on a flag variety --- a ``generic torus orbit
closure'', see   \cite{kn:dab}  and 
Theorem  \ref{thmdabrow}.

In this paper we  generalize Voskresenski\u{\i} and Klyachko
results: we describe  all pairs $(R,L)$ such that the associated toric
variety $X_{R,L}$ is   $\QQ$-Gorenstein Fano and which varieties
among them are  Gorenstein Fano (see Definition
  \ref{defi:QGFGFF}).  
The smoothness condition (as in \cite{kn:voskly}) imposes restrictions of the
combinatorics, in such a way that (in the case of irreducible root
systems) smooth Fano varieties are obtained only for root systems of type
$A_n$ (two infinite series), $C_n$ (one additional infinite series)
and $G_2$ (one exceptional case). By relaxing  the
smoothness constraint, we obtain 
$\QQ$-Gorenstein Fano generic closures  for all types of
irreducible root systems except $E_7$ and $E_8$ (see Section
\ref{section:table}). More precisely, we 
exhibit  twelve
infinite series plus five exceptional cases of $\QQ$-Gorenstein Fano
varieties. 

From a combinatorial perspective,
the exhibition of toric Gorenstein Fano varieties is  interesting because 
each of these varieties is naturally associated with a couple of dual
reflexive polytopes  (see
Definition \ref{defn:latpol}) --- this duality allowed Batyrev 
to give a rigorous construction of mirror
symmetry in the toric context, see \cite{kn:BatyMiroir}. 
Applying in our setting this well known correspondence, we produce the
list of  dual reflexive
polytopes associated  to each Gorenstein Fano toric
  variety $X_{R,L}$.   As a minor
by-product, we describe the list of  root 
polytopes  which are reflexive.

We briefly describe now the content of this paper.
In Section \ref{sec:prelim} we establish the basic notations and
we present some known results on Gorenstein Fano varieties and  closures of
generic orbits.
In Section \ref{sec:Combinatorics} we study some combinatorial
properties of the cone $\sigma_{R,L}$ and we  characterize the 
$\QQ$-Gorenstein Fano generic closures in terms of the combinatorics
of $\sigma_{R,L}$ (see Theorem
\ref{th:FinalCriteria}).
In  Section \ref{section:table} we present our main result (Theorem
\ref{th:Classification}), namely a classification of all
$\QQ$-Gorenstein Fano,  Gorenstein Fano and smooth Fano generic closures in
terms of their defining set of roots $L\subset I$ (see Definition \ref{defi:QGFGFF}). The proof of this
classification result 
relies heavily in the criterion established in Section  \ref{sec:Combinatorics}.
In Section \ref{sec:Reflexive} we exhibit the reflexive polytopes
associated  to the Gorenstein Fano generic closures;  as an application of our
combinatorial description, we classify the reflexive  root polytopes.  

In \cite{kn:code},
the interested reader can find a SAGE (see \cite{sage}) package 
that allows to perform explicit calculations for the cone
$\sigma_{R,L}$ --- we use    GAP3   (version 
maintained by Jean Michel, \cite{kn:gap3}) in order  to use the package Chevie
(see \cite{kn:chevie} and  \cite{kn:chevieMichel}).

\section{Preliminaries}
\label{sec:prelim}

\subsection{$\QQ$-Gorenstein Fano toric varieties}\label{sec:GFtoric}\label{sec:toric}\ %

 All along this work, by a  toric variety we mean  a normal toric
  variety over an algebraically closed field $\Bbbk$; our general
  reference for toric varieties is \cite{kn:cox}.

If  $T$ is an algebraic torus, we denote by  $\Lambda$  the characters
group of $T$ and by $\Lambda^\vee$ the
$\ZZ$-dual of $\Lambda$.  We denote  by
$\Lambda_\QQ$ 
(resp. $\Lambda^\vee_\QQ$) the $\QQ$-vector space $\QQ\otimes_{\mathbb
	Z} \Lambda$ (resp. $\QQ\otimes_{\mathbb Z} \Lambda^\vee$), and if
$(u,v)\in \Lambda_\QQ\times\Lambda^\vee_\QQ$ then
$\Dual{u}{v}=v(u)\in\QQ$ is the natural pairing of $u$ and $v$.  

If $X$ is a subset of a finite dimensional $\QQ$-vector space $V$, we
denote by $\Conv(X)$ the convex hull of $X$, by $\langle X\rangle$ the vector space
generated by $X$, and by $\langle X\rangle_{\operatorname{aff}}$ the affine space
generated by $X$; $\dim_\aff(X)$ denotes the affine dimension of $\langle X\rangle_{\operatorname{aff}}$. We denote  by $\QQ^+X$ the positive cone generated by
$X$ (with the origin as  vertex); the   ``dual cone''  of
$X$ is defined as 
\[
X^\vee:=\bigl\{\varphi\in V^\vee\mathrel{:}\forall x\in
X,\,\Dual{x}{\varphi} \geq 0\bigr\}. 
\]

% \begin{definition}
	If $\Sigma$ is a fan in $\Lambda^\vee_\QQ$ (see
        \cite[Definition 3.1.2]{kn:cox}), then  $\Sigma(r)$ is the set
        of $r$-dimensional cones in
	$\Sigma$. For each $\rho\in\Sigma(1)$, 
        $u_\rho$ is  the 	\emph{primitive element} of 
	the  monoid $\rho\cap\Lambda^\vee$.  
	  The set of primitive 	elements of $\sigma\in \Sigma$ is
          denoted by
	\[\Prim{\sigma}=\bigl\{ u_\rho\mathrel{:}
	\rho\in\Sigma(1)\mbox{ and } \rho\subset\sigma\bigr\}. 
      \]
      
	We set
        $\Prim{\Sigma}=\bigcup_{\sigma\in\Sigma}\Prim{\sigma}$.
        
	A fan  $\Sigma$ has associated a toric variety  that we denote as
	$X_\Sigma$.

      % If $\sigma$ is a given rational strictly convex cone, then the
      %   set of its faces define  a fan that we will denote also by $\sigma$.
% \end{definition}

% \begin{remark}

	% Recall that if $\Sigma$ is a fan in $\Lambda^\vee_{\QQ}$, then:
	
	% \noindent (1) the toric
	% variety $X_\Sigma$ 
	% is smooth if and only if for each $\sigma\in\Sigma$, $\sigma$ is
	% generated by a finite set $S_\sigma\subset  \Lambda^\vee$ which can be completed to a
	% basis of $\Lambda^\vee$;
	
	% \noindent (2) the variety $X_\Sigma$ is  complete if and only if
	% $\Sigma$ is complete {\it i.e.} 	$\bigcup_{\sigma 
	% 	\in \Sigma}\sigma=(\Lambda^\vee)_\QQ$;
	
	% \noindent (3)  the cones in  $ \Sigma(r)$ are in bijection
	% with the  irreducible closed $T$-varieties of $X_{\Sigma}$ of
	% co-dimension $r$. In 
	% particular, the $T$--stable Weil divisors are  in one-to-one
	% correspondence with the $\ZZ$--linear combinations $\sum_{\rho\in\Sigma(1)}a_\rho
	% D_\rho$, where $a_\rho\in \ZZ$ and $D_\rho$ is the $T$-stable divisor
	% associated to the cone $\rho\in \Sigma(1)$.
% \end{remark}

% \begin{definition} 
	Recall that if  $\sigma \subset \Lambda_\QQ^\vee$ is a
	polyhedral strictly convex
	cone, then the \emph{relative interior} of $\sigma$, denoted by
	$\mathring{\sigma}$,  is the complement in $\sigma$ of the union of the facets of $\sigma$. 
% \end{definition}

 \begin{definition}
 	\label{defn:latpol}
	If $\Lambda$ is a lattice, a \emph{lattice polytope} is the convex
	hull in $\Lambda_\QQ$ of a finite subset $X \subset \Lambda$.
        
 Assume now that   $\mathcal P\subset\Lambda_\QQ$ is a lattice
        polytope of maximal rank         containing the  origin in its
        strict interior. If  $\mathcal Q$ is a proper  facet of
        $\mathcal P$, the \emph{interior normal} 
	of $\mathcal Q$, denoted by  $\varphi_{\mathcal Q}$, is the unique
	element of  $(\Lambda_\QQ)^\vee$ such that  $\varphi_{\mathcal
          Q}(Q)=-1$ 	 and $\varphi_{\mathcal Q}(\mathcal
        P\setminus\mathcal Q)>-1$; the  \emph{exterior normal} is defined as
        $-\varphi_{\mathcal Q}$. The  convex hull of the set of the 
	interior normals of $\mathcal
        P$ is called the \emph{dual
        polytope of  $\mathcal   	P$} and is denoted by
	$
	\mathcal P^\circ=\{u\in\Lambda^\vee_\QQ\mathrel{ :
        }\Dual{v}{u}\geq -1\ \forall v\in\mathcal P \}$. 
	
	It is clear that  $\mathcal P^\circ$ is a polytope and that $(\mathcal
	P^\circ)^\circ=\mathcal P$; we say that 
	$\mathcal P$ is a \emph{reflexive polytope} if $\mathcal P^\circ$ is
	a lattice polytope.
 \end{definition}

 % \begin{remark}\label{rq:DescriptionDual}
 %        If  $\mathcal P$ is a lattice polytope  of
 %        maximal rank  containing the origin in its interior, then 
 %         \marginpar{\color{blue}formule elim}
 %        % that is: 
 %        % \[
 %        % \mathcal P^\circ\=\Conv \bigl(\{\varphi_{\mathcal Q}\mathrel{ : } \mathcal
 %        % Q\mbox{ is a facet of }\mathcal P\}\bigr).
 %        % \]
        
 % \end{remark}

\begin{definition}\label{defi:QGFGFF}
  Let $X$ be a normal variety and  denote by
  % $K_X$ the canonical
%	divisor
        %--- that is the divisor associated to the canonical sheaf ---
        % and by
  $-K_X$ the 
  anti-canonical divisor.
  %(associated to the dual of the canonical sheaf).	
	 We say that  $X$ is  $\QQ$-Gorenstein Fano  if $-K_X$ is
        an ample $\QQ$-Cartier divisor;  if moreover
          $-K_X$ is an ample Cartier divisor,  we say that $X$ is
          Gorenstein Fano. If $X$ is a smooth Gorenstein Fano
        variety, we say that $X$ is smooth Fano. 
	
	If $X$ is $\QQ$-Gorenstein Fano, the \emph{Gorenstein index} of $X$, denoted by $j_X$,  is
	the  smallest positive integer $j$ such that $jK_X$ is
        Cartier
        --- thus, a Gorenstein Fano variety is a $\QQ$-Gorenstein
        Fano variety  of Gorenstein index $1$.
	
%	A smooth variety $X$ is  Fano  if  $-K_X$ is an ample divisor.

	A fan  $\Sigma\subset\Lambda^\vee_{\QQ}$ is $\QQ$-Gorenstein Fano,
	Gorenstein Fano or smooth Fano if the associated toric variety
        $X_\Sigma$ has the corresponding property.    
\end{definition}

% \begin{remark}
% 	Recall that if   $\Sigma\subset\Lambda^\vee_{\QQ}$ is a fan then
% 	$-K_{X_\Sigma}=\sum_{\rho\in\Sigma(1)}D_\rho$. 
% \end{remark}

The following equivalences are well known (see for example
\cite[Theorem 4.2.8 and Lemma
6.1.13]{kn:cox}): 

\begin{proposition}
	\label{prop:fanogeneral}
	Let $\Sigma$ be a complete fan in $\Lambda^\vee_\QQ$.  The following assertions are
	equivalent:
	\begin{enumerate}[(1)]
		\item  $X_\Sigma$ is  a $\QQ$-Gorenstein
                  Fano toric variety;
		
		\item $\bigl\{
		\Conv\bigl(\Prim{\sigma}\bigr) \mathrel{:} \sigma\in\Sigma(s)\,,\ s=1,\dots,
		n\bigr\}$  is the set of  proper faces of the lattice polytope 
		$\Conv\bigl(\Prim{\Sigma}\bigr)$;
		
		\item for every cone $\sigma\in \Sigma(n)$, the
                   polytope
		$\Conv\bigl(\Prim{\sigma}\bigr)$ is $(n-1)$-dimensional; let
		$\varphi_\sigma\in\Lambda_\QQ$  be such that 
		$\Dual{\varphi_\sigma}{v}=-1$ for $v\in \Prim{\sigma}$. Then $\Dual{\varphi_\sigma}{w} >
		-1$ for every $w\in \Prim{\Sigma}\setminus\Prim{\sigma}$. 
	\end{enumerate}	
	Moreover, if  $X_\Sigma$ is a $\QQ$-Gorenstein Fano
        toric variety, then
	\[
	j_X=\min\bigl\{ j\in \NN^*\mathrel{:}  \forall u\in\Lambda^\vee\,,\ \forall\sigma\in
	\Sigma(n)\,,\Dual{j\varphi_\sigma}{u}
	\in\ZZ \bigr\}.
	\]

	In particular, if $X_\Sigma$ is Gorenstein Fano then
        $\Conv\bigl(\Prim{\Sigma}\bigr)$ is a reflexive lattice polytope.\qed
\end{proposition}

%\section{Toric variety associated to Root systems}\label{sec:RootSystems} %

\subsection{Fans defined by root systems and generic orbits}\ %
\label{subsec:RootSystems}

 In this section we  establish our notations on fans
  defined by root systems, and formulate  Dabrowski's results
  accordingly, associating to a generic orbit the corresponding combinatorial data.
When dealing with root systems, we follow
  Bourbaki's notations  (see \cite{kn:bourbaki,kn:bourbakien}).

\begin{notation}
In what follows,  $R$ designs a root system of rank $n$ and
$\Lambda_R$ its associated \emph{root lattice}. We denote by
$R^+$ a chosen set of  positive  
roots;   $S=\{\alpha_i : i\in I\}$ is   the set of the \emph{simple
  roots} associated to $R^+$.
%  --- this set of simple roots is a basis of  the lattice $\Lambda_R$.
We denote by $\Lambda_P$ the lattice of weights and by
$\{\omega_i : i\in I\}$ the set of  \emph{fundamental  
  weights}  associated to $S$.
% ---  the set of fundamental weights  is a basis of the weight lattice $\Lambda_P$.

If $\alpha\in R$, we denote by $s_\alpha:(\Lambda_P)_\QQ\to
(\Lambda_P)_\QQ$ the associated reflection, and by $W$ the Weyl group
 generated by the reflections associated to $R$.
 Recall that $W$  acts on $(\Lambda_P)_\QQ$ with the Weyl chamber $\mathcal
D= \QQ^+\{\omega_i \mathrel{:} i\in I\}$ as
fundamental domain.

The root system dual to $R$ is denoted by $R^\vee$. Recall that the  \emph{simple
  co-roots} $\{\alpha^\vee_i \mathrel{:} i\in I\}\subset 
(\Lambda_P)_\QQ^\vee$ and  the \emph{fundamental co-weights}
$\{\omega^\vee_i\mathrel{:} i\in I\}\subset (\Lambda_R)_\QQ^\vee$
are  such that $\Dual{\alpha_i^\vee}{\omega_j}=\delta_{ij}$, 
$\Dual{\alpha_i}{\omega_j^\vee}=\delta_{ij}$  for all
$(i,j)\in I^2$.  Also, the  reflections $s_{\alpha_i^\vee}$ induce an action
of $W$ on $(\Lambda_P^\vee)_\QQ$, with the dominant Weyl Chamber of
$R^\vee$ (denoted by   $\mathcal
D^\vee$) as fundamental domain.
	% If the root system is \emph{simply laced}, the  scalar
	% products $\Scalar{\cdot}{\cdot}$ are chosen so that
	%         $\Lambda_R^\vee$ is identified with 	$\Lambda_P$.

        % . In this case, we  identify $\Lambda_R^\vee$ with
	% $\Lambda_P$ by means of the scalar product --- recall that under
	% this identification $R^\vee=R$. 

If $L\subset I$, we will abuse notations and identify
  $S_L=\{\alpha_i\mathrel{:} i\in
  L\}$ with $L$.
We  denote by  $W_L$ 
	the subgroup of $W$ generated by the corresponding reflections 
	$\{s_{\alpha_i^\vee}\mathrel{:} i\in L\}$.

We denote by   $R=\prod_{k=1}^r R_k$ the decomposition of the root
system $R$ in irreducible root systems;  the set
of simple roots of $R_k$ is denoted by 
$S(k)\subset S$ and we denote by $I(k)\subset I$ the corresponding subset of indexes.

 We choose $W$-invariant scalar products in $(\Lambda_P)_\QQ$ and 
$(\Lambda_P^\vee)_\QQ$; these
scalar products will be denoted by $\Scalar{\cdot}{\cdot}$ in both cases. 

      \end{notation}

\begin{definition}
	With the previous notations, if $\lambda=\sum_{i=1}^n a_i\omega_i\in
	\Lambda_P$ is a  weight (resp.  $\lambda^\vee=\sum_{i=1}^n a_i\omega^\vee_i\in
	\Lambda_R^\vee$ is a co-weight), we define the \emph{support} of
	$\lambda$ (resp. $\lambda^\vee$) as the set
        $I_\lambda=\{\alpha_i\mathrel{:} a_i\neq 0\}\subset S$.
	\end{definition}

From now on, $G$ is a semi-simple group over  $\Bbbk$ and  $T\subset G$ a maximal
torus, such that  $R$ is the root system associated to the couple
$(G,T)$; we denote by $B\subset G$ the  Borel subgroup associated to
$R^+$. To each subset
$L\subset I$, 
we associate the parabolic subgroup $P_L$ containing the opposite
Borel subgroup  $B^-$ and  such that  the Weyl group of $P_L$  is
equal to $W_L$.

Recall that if $\lambda\in \Lambda_{P}$ is  a dominant weight with  support
$I_\lambda$ contained in $L^c=I\setminus L$, then $\lambda$ can be extended to $P_L$. We denote by 
$V(\lambda)$ the \emph{Weyl $G$--module} associated to $\lambda$.
% --- recall that 
% $V(\lambda)=\bigl\{ f\in \Bbbk [G] \mathrel{:} 
% f(xy)=\lambda^{-1}(y)f(x)\ \forall x\in G, y\in P \bigr\}$.

\begin{definition}\label{def:FanRoot}
	\label{def:SigmaRL}
	Let $L\subsetneq I$ be a proper subset of roots.  
        Following \cite{kn:voskly}, we define the \emph{cone
		associated to $L$} as  
	\[
	\sigma_{R,L}=\bigcup_{w\in W_L}w\mathcal D^\vee\subset (\Lambda_R)_\QQ^\vee.
	\]

        If $\sigma_{R,L}$ is strictly convex, we consider the
 complete  fan having as maximal  
	cones  the translates $w\cdot
	\sigma_{R,L}$, where $w\in W$; we denote this fan 
           by 	$\Sigma_{R,L}\subset
           (\Lambda_R)^\vee_\QQ$.

           We define    $-\Sigma_{R,L}:=\bigl\{ -\sigma\mathrel{:}
	\sigma\in\Sigma_{R,L}\bigr\}$; the corresponding toric variety
        is denoted by $X_{R,L}:=   X_{-\Sigma_{R,L}}$.

%           denote by 	$\Sigma_{R,L}\subset
%           (\Lambda_R)^\vee_\QQ$   the complete  fan having as maximal  
% 	cones  the translates $w\cdot
% 	\sigma_{R,L}$, where $w\in W$. Notice  that if $W^L$ is a
%         set-theoretical section of $W/W_L$, then 
% $	\Sigma_{R,L}(n)=\bigl\{w\cdot \sigma_{R,L}\mathrel{:} w\in
% 	W^L\bigr\}$.
	% We denote by    $-\Sigma_{R,L}=\bigl\{ -\sigma\mathrel{:}
	% \sigma\in\Sigma_{R,L}\bigr\}$  and define  $X_{R,L}=   X_{-\Sigma_{R,L}}$.

	% If $\sigma_{R,L}$ is a strictly convex cone,
        %   then         all the translates
	% of $\sigma_{R,L}$ by the elements of $W$ are so, and $w\cdot
	% \sigma_{R,L}\cap w'\cdot \sigma_{R,L}$ is a common face for all
	% $w,w'\in W$;  we denote by 
	% $\Sigma_{R,L}\subset (\Lambda_R)^\vee_\QQ$   the fan having as maximal 
	% cones these translates:

      % \[
      %   \Sigma_{R,L}(n)==\bigl\{w\cdot \sigma_{R,L}\mathrel{:} w\in W^L\bigr\}, 
      %   \]
      %   where $W^L\subset W$ is a set-theoretical section of $W/W_L$.

	% Notice that $\Sigma_{R,L}$  is a complete fan
	% and that $-\Sigma_{R,L}=\bigl\{ -\sigma\mathrel{:}
	% \sigma\in\Sigma_{R,L}\bigr\}$  is also a complete fan.  We denote $X_{R,L}=   X_{-\Sigma_{R,L}}$.
	
\end{definition}

\begin{remark}
	\label{rem:prodfan}
\noindent \emph{(1)}	
		The geometric meaning of
		 the use of the
		co-weight lattice and the fan $-\Sigma_{R,L}$ in the
                definition of $X_{R,L}$  
		(e.g. instead of the weight lattice and $\Sigma_{R,L}$),  will become
		evident in the next section (see Theorem 
		\ref{thmdabrow} and Remark \ref{rem:probdrab}).\\

\noindent \emph{(2)}	 Let $R=\prod_{k=1}^rR_k$ be a decomposition of the root
		system $R$ in irreducible root systems. Then  $\mathcal
		D^\vee=\prod_{k=1}^r \mathcal D_k^\vee$, 
		where $\mathcal D_k^\vee\subset (\Lambda_k)^\vee_{\QQ}$ is the
		dominant Weyl chamber of $R_k^\vee$ --- here $\Lambda_k$ denotes the
		lattice generated by $R_k$.
			In particular, $\sigma_{R,L}=\prod_{k=1}^r\sigma_{R_k, L\cap I_k}$, $\Sigma_{R,L}=\prod_{k=1}^r\Sigma_{R_k,L\cap I_k}$,
and   $X_{R,L}\cong \prod_{k=1}^rX_{R_k,L\cap I_k}$.
	
\end{remark}

 Dabrowski proved in \cite{kn:dab} that the toric varieties
 $X_{R,L}$ can be realized as closures of ``generic'' $T$-orbits in $G/P_L$. We
briefly  recall  his construction, filling some minor gaps in the proofs
 presented in \emph{op.cit.}   for the
sake of  completeness.

\begin{definition}[{see \cite[{\textsection 1}]{kn:dab}}]
	Let $L\subset I$ be a  subset of roots,  $\Pi_\lambda=\bigl\{\mu \in\Lambda_P \mathrel{:} V(\lambda)_\mu\neq 0\bigr\}$ the set of 
	$T$-weights of $V(\lambda)$,  and $\mathcal A_\lambda$ be the list of the $T$-weights counted 
	with multiplicity. A \emph{set of Pl\"ucker coordinates}
        $\{f_\mu\mathrel{:} \mu\in \mathcal  
	A_\lambda\}$  is a choice of a basis of $T$-semi-invariants functions
	$f_\mu\in V(-\lambda)_\mu$.  

        If $x=uP\in G/P$, we consider
	$
	\Pi_\lambda(x):=\bigl\{\mu \in\Pi_\lambda\mathrel{:} f_\mu(x)\neq 0 \text{ for some } f_\mu\bigr\}$. 
	% It is easy to see that $\Pi_\lambda(x)$ does not depend on the choice of the Pl\" ucker 
	% coordinates. Moreover, $\lambda -w\Pi_\lambda(x)\subset S^L\subset \Lambda_R$, for every
	% $w\in W$.
		We say that the $T$-orbit $T\cdot x$ is \emph{generic}
                in the sense of Dabrowski if 
	$W\cdot \lambda\subset \Pi_\lambda(x)$ and the set $\lambda
        -w\Pi_\lambda(x)$ generates $S^L$ as a sub-monoid --- notice 
        that if all the Pl\"ucker coordinates of $x$ are non zero,
        then $T\cdot x$ is generic.
\end{definition}

\begin{proposition}\label{prop:dualcone} \label{prop:sigma.convex}
	Let $(R^L)^+$ be the set of positive roots which are not sum
	of simple roots  in $L$, and $S^{L}$ be the
	sub-monoid generated by $(R^L)^+$.  
	Then $\sigma_{R,L}$ is the dual cone of the convex cone
        generated by  $S^{L}$. 
In particular,   $\sigma_{R,L}$ is  a strictly convex  cone if and only if 
	$L\cap I(k)\neq I(k)$ for all  $k=1,2,\ldots,r$.

      \end{proposition}
\begin{proof}
	Since $
	\QQ^+{(S^{L})}=\QQ^+{\bigr((R^{L})^+\bigl)}$, it follows that
	\[
	\bigl(\QQ^+(S^L)\bigr)^{\vee}= \bigl(\QQ^+\bigl((R^L)^+\bigr)\bigr)^{\vee} 
        =\bigcap_{\beta\in (R^L)^+} \bigl\{\chi^\vee\in(\Lambda^\vee_R)_\QQ\mathrel{:}\Dual{\beta}{\chi^\vee}\geq 0\bigr\}.\]
	
	It is clear that  $(R^L)^+$ is $W_L$-stable; hence, 
	$\QQ^+\bigl((R^L)^+\bigr)^{\vee}$  is also
        $W_L$-stable  and, as it contains 
	the dominant chamber $\mathcal D^\vee$, we have the inclusion
	$\sigma_{R,L}\subset \QQ^+\bigl((R^L)^+\bigr)^{\vee}$.  In order to prove the
	equality, it suffices to prove that 
	$s_{\alpha_i}\cdot \sigma_{R,L}$ is not contained in
	$\QQ^+\bigl((R^L)^+\bigr)^{\vee}$ for all $i\notin L$. But if $\langle
	\alpha_i,\chi^\vee\rangle > 0$, then
	$
	\bigl\langle   \alpha_i,s_{\alpha_i}(\chi^\vee)\bigr\rangle= - \langle
	\alpha_i,\chi^\vee\rangle< 0$.

	In order to prove the converse, in view of Remark \ref{rem:prodfan}, we can assume that $R$ is
	irreducible. In this case, it is clear that
	$\sigma_{R,I}=(\Lambda_R)^\vee_\QQ$. If $L\neq I$, then
	it is well-known that  $(R^L)^+$ generates a space of maximal
	dimension and the result follows.
      \end{proof}

 We are in condition now to state Dabrowski's main result concerning
 the generic $T$-orbits in $G/P_L$ and their closures.

\begin{theorem}[{\cite[Theorem 3.2]{kn:dab}}]\label{thmdabrow}
	If $L\subset I$ is such that $L\cap I(k)\neq I(k)$ for all $k$
        (see Proposition \ref{prop:sigma.convex}) and
          $T\cdot x\subset G/P_L$ is a generic orbit, then $\overline{T\cdot
		x}$ is a toric variety isomorphic to the toric variety $ X_{R,L}$.\qed

        % In particular, generic orbits exist.

	% If $T\cdot x $ is a
	% generic orbit, then the $T$-orbit closure   $\overline{T\cdot
	% 	x}\subset G/P$ is a toric variety. This toric variety is isomorphic
	% to the variety $X_{R,L}$ defined by the fan $-\Sigma_{R,L}$. 
	
\end{theorem}

\begin{remark}
	\label{rem:probdrab}
%	We keep the previous notations.
	
	\noindent (1)   Since the torus $T$ does not
	act  effectively on $G/P$, in general the generic orbit associate to
	$L$ has strictly lower dimension than $T$ --- that is $T\not\cong
	T\cdot x$.
	
	However, it is easy to show that in the hypothesis of Proposition
	\ref{prop:sigma.convex}  we have that  $T_x=Z(G)$. Since 
	$\Lambda_R$ is the lattice of characters of $T/Z(G)$, the
	description of $\overline{T\cdot x}$ as a $T/Z(G)$-toric variety is
	given by a complete fan  in  the
	space $(\Lambda_R^\vee)_\QQ$.
	
	\noindent (2) When the root system is simply laced, one has some
	leeway for the combinatorial description of a 
	closure of a generic orbit as a toric variety, since
        $\Lambda_R^\vee=\Lambda_P$ under the  identification
        of $R$ and $R^\vee$.
        % (see Remark \ref{rq:SimplyLaced}).
	This discretionality
	appears in the literature, sometimes by omission, e.g. in \cite{kn:dab}
	the lattice is not explicitly mentioned; however, the reader should be
	aware that if the root system is not simply laced, then $\Lambda_P$ and 
	$\Lambda_R^\vee$ are distinct lattices. The  distinction between  $\Lambda_P$ and 
	$\Lambda_R^\vee$ must be taken into account in order to give a correct combinatorial  description of the fan associated   to adherence of a generic orbit, see for example  \cite{kn:Klyach1}.
	% Notice that in the case of $L=\emptyset$, this right description of the fan associated to the adherence to the generic orbits is given in \cite{kn:Klyach1}.
\end{remark}

\section{A criteria for $X_{R,L}$ to be $\QQ$-Gorenstein Fano}\label{sec:Combinatorics}

In this section  we characterize
when the  closure of a  generic orbit  is a Gorenstein
  Fano toric variety in 
terms  of the combinatorial properties of the core of the associated
cone  $\sigma_{R,L}$. From now on we assume  that $L\subset I$ is such that $\sigma_{R,L}$ is a strictly
convex cone.

%  We gather here various properties of the  $\sigma_{R,L}$ which
% are useful for the main result of this paper.   
% First we describe the relationship  between $\Sigma_{R,L}$ and the
% Weyl polytope associated to an  associated  dominant weight. Then we
% introduce the \emph{core} of $\sigma_{R,L}$ and present some basic
% properties of the core --- that will simplify our characterization of
% the generic orbits having Gorenstein Fano closure.

\subsection{Various combinatorial properties of $\sigma_{R,L}$}\ %
\begin{definition}\label{def:WeylPolytope}
	If  $\lambda\in(\Lambda_R)_\QQ$, the 
	 \emph{Weyl polytope associated to $\lambda$} is defined as
	$\mathcal{WP}(\lambda)=\Conv(W\cdot \lambda)\subset (\Lambda_P)_\QQ$. 
The set of facets of $\mathcal{WP}(\lambda)$  containing $\lambda$ is
denoted by $C_\lambda(n-1)$.
      \end{definition}

In \cite{kn:dab}, Dabrowski showed that  $\Sigma_{R,L}$
  can be obtained as the fan dual of the Weyl polytope of a dominant weight
   with support $L^c=I\setminus L$:

\begin{proposition}\label{prop:WeylPolytope}
	Let $\lambda$ be a dominant weight with support $L^c$, and assume that
	$\sigma_{R,L}$ is strictly convex. Then
	\[
          \Lambda_R\cap\bigl(\lambda-\mathcal{WP}(\lambda)\bigr)=S^L\cap\bigl(\lambda-\mathcal{WP}(\lambda)\bigr).
          \]
	  
 Moreover, the fan $\Sigma_{R,L}$ is dual to the polytope
  $\mathcal{WP}(\lambda)$. \qed 
\end{proposition}

The cones  $\sigma_{R,L}$ being stable by the action of $W_L$, we can
use this action in order to describe their combinatorics and geometry
as follows.

\begin{proposition}\label{prop:faceSigmaLambda}
%  Let  $L\subset I$
%  and consider the associated cone $\sigma_{R,L}$.
%   Then
%   the group $W_L$ acts  on the set of faces 
% 	$\sigma_{R,L}(r)$, and any face $\mathcal F\in\sigma_{R,L}(r)$ is the union of the translates by
% 	$W_L$ of faces of $\mathcal  D^\vee(r)$ contained in $\mathcal F$.
% Let  $L\subset I$
%   and consider the associated cone $\sigma_{R,L}$.
	
	% \begin{enumerate}[(1)]
	% 	\item $\mathcal F=\bigcup_{i=1}^{s} w_i\cdot\gamma_i$. 
		
	% 	\item For all  $i, j=1,\dots,s$, {\color{blue}$i\neq j$,} $  w_i\cdot \gamma_i\cap
	% 	w_j\cdot \gamma_j$ is a common proper face of $w_i\cdot \gamma_i$ and $w_j\cdot \gamma_j$.
		
	% 	% \item If $w\in W_L$ and $\gamma\in \mathcal D^\vee(r)$, are such
	% 	% $w\cdot \gamma \subset \mathcal F$, then $w=w_i$ and $\gamma=\gamma_i$ for
	% 	% some $i\in\{1,\dots ,s\}$. In  particular, if $w_i=w_j$ then
	% 	% $\gamma_i=\gamma_j$.

        %       \end{enumerate}

 	% In particular,
	% \[
	% \Prim{\sigma_{R,L}}=W_L\cdot \bigl(\Prim{\sigma_{R,L}}\cap
	% \{\omega^\vee_1,\dots, \omega_n^\vee\}\bigr).
	% \]

    Let  $L\subset I$ and consider $\mathcal F\in\sigma_{R,L}(r)$. Then 
	there exist  unique pairs $(\gamma_i,w_i)\in
        \mathcal D^\vee(r)\times W_L$, $i=1,\dots,s$,  such that
        $\mathcal F=\bigcup_{i=1}^{s} w_i\cdot\gamma_i$. 
 Moreover,
        \begin{enumerate}[(1)]
        \item for all  $i, j=1,\dots,s$, $i\neq j$, $
w_i\cdot \gamma_i\cap 
		w_j\cdot \gamma_j$ is a common proper face of
                $w_i\cdot \gamma_i$ and $w_j\cdot \gamma_j$;

\item  if $w_i=w_j$ then
		 $\gamma_i=\gamma_j$.

              \end{enumerate}

	In particular, $
	\Prim{\sigma_{R,L}}=W_L\cdot \bigl(\Prim{\sigma_{R,L}}\cap
	\{\omega^\vee_1,\dots, \omega_n^\vee\}\bigr)$.

\end{proposition}

\begin{proof}
	By construction, $\sigma_{R,L}$ is stable under the (linear) action of $W_L$;
	hence  $\sigma_{R,L}(r)$ is stable under the action of
	$W_L$. Since $\mathcal D^\vee$ is a  fundamental domain
	of the action of $W$
        %and that the action is linear with finite orbits,
	it follows that
	there exist $w_1,\dots w_s\in W_L$ and $\gamma_1,\dots,
	\gamma_s\in\mathcal D^\vee(r)$ such that  $\mathcal F=\cup_{i=1}^sw_i\cdot
	\gamma_i$.  In particular,   it
            follows  that   $\mathcal D^\vee$ is also  a 
	fundamental domain for action  of $W_L$ on $\sigma_{R,L}$.
	
	Let $w\in W_L$, $\gamma\in\mathcal D^\vee$ be such that $w\cdot
	\gamma\subset 
	\mathcal F$. Since the affine dimension of 
          $w\cdot \gamma$ and $w_i\cdot
	\gamma_i$, $i=1,\dots ,s$,    is $r$, it follows
	that there exists $i\in\{1,\dots,s\}$ such that $\dim \langle w_i\cdot
	\gamma_i\cap 
	w\cdot \gamma\rangle_\aff = r$.   Since $\mathcal D^\vee$ is a 
	fundamental domain for the action of $W_L$ on $\sigma_{R,L}$,
	it follows that  $w_i\cdot
	\gamma_i=w\cdot \gamma$ and therefore $w_i=w$ and
	$\gamma_i=\gamma$.
	
	It is now easy to prove assertion (1). In order
          to prove (2), observe that if $\gamma_i\neq
        \gamma_j$, then $\dim \langle \gamma_i\cup
        \gamma_j\rangle_\aff>r$ and it follows that $w\cdot (\gamma_i\cup
        \gamma_j)$ is not included in $\sigma_{R,L}(r)$ for any $w\in W_L$.

        Finally, recall that $\Prim{\sigma_{R,L}}$ consists of the primitive
	elements of the rays in  $\sigma_{R,L}(1)$ and that the lattice of
	co-roots is stable under de action of $W$.
\end{proof}

\begin{notation}\label{nota:JL}
	We denote $J_L=\Prim{\sigma_{R,L}}\cap
	\{\omega^\vee_1,\dots, \omega_n^\vee\}$ --- by Proposition
	\ref{prop:faceSigmaLambda} above, $J_L$ is a fundamental domain for the
	action of $W_L$ on $\Prim{\sigma_{R,L}}$. In order to simplify the
	notations, we will often also denote the set of indexes $\bigl\{
	i\mathrel{:} \omega_i^\vee\in J_L\bigr\}$ by $J_L$. 
\end{notation}

Once we have described $\Prim{\sigma_{R,L}}$  as a set of
$W_L$--orbits, we are in condition to give a simple description of
its affine support space.

\begin{proposition}\label{prop:DescriptionNormal}
	If  $L\subset I$ then
	 \[
	\begin{split}
 	\bigl\langle \Prim{\sigma_{R,L}}\bigr\rangle_\aff= & \ \scriptstyle \omega_k^\vee +
	\bigl\langle \bigl(\bigcup_{j\in J_L} W_L\cdot (\omega_j^\vee)
	-\omega_j^\vee \bigr)\, \cup\,   \{ \omega_i^\vee -\omega_j^\vee\mathrel{:}\,  i,j\in J_L \}\bigr\rangle\\
 \ 	=& \  \scriptstyle  \omega_k^\vee + \bigl\langle \{ \alpha^\vee_i\mathrel{:}\, i\in L \}\cup  \{
	\omega_i^\vee -\omega_k^\vee \mathrel{:}\, i\in J_L\}\bigr\rangle
	\end{split}
	\]
	for any $\omega^\vee_k\in J_L$.
\end{proposition}

\begin{proof} 
	By Proposition \ref{prop:faceSigmaLambda},
	$\Prim{\sigma_{R,L}}=W_L\cdot J_L$.  If
	$\omega^\vee_i,\omega_j^\vee\in J_L$ and $w,w'\in W_L$, then 
	$
	w\cdot(-\omega^\vee_i)-w'\cdot(-\omega^\vee_j)= w\cdot(-\omega^\vee_i)-\omega_i +\omega^\vee_i-\omega^\vee_j+\omega^\vee_j-  w'\cdot(-\omega^\vee_j)$,
	and the first  equality follows. 
	
	For the second equality, let $w=s_\ell\cdots s_1\in W_L$,
	with $s_j\in \{ s_{\alpha_i^\vee}\mathrel{:} i \in L\}$. Then
	$w\cdot(-\omega^\vee_t)-\omega^\vee_t\in \langle \alpha^\vee_i\mathrel{:} i\in
	L\rangle_\QQ$ for all $\omega_t^\vee\in J_L$  and the
	inclusion $\subset$ follows.
	
	In order to prove the remaining inclusion, let $i\in L $. If
	$s_{\alpha_i^\vee}\cdot \nu=\nu$ for all $\nu\in 
	\Prim{\sigma_{R,L}}$, then $s_{\alpha_i^\vee}$ acts trivially, since 
	$\sigma_{R,L}$ is of maximal dimension; this is a contradiction. It
	follows that there exists $\nu\in \Prim{\sigma_{R,L}}$ such that
	$s_{\alpha_i^\vee}\cdot \nu\neq \nu$, and therefore $\alpha^\vee_i\in \bigl\langle \Prim{\sigma_{R,L}}\bigr\rangle_\aff-\omega_k^\vee$.
\end{proof}

 Since by Proposition \ref{prop:DescriptionNormal}  the set
 	$J_L$ determines $\Sigma_{R,L}$ and therefore $X_{R,L}$, we
        proceed to calculate $J_L$, by 
       translating results by Khare to our context (see \cite[Definition
       3.1 and Theorem C]{kn:khare}).

        \begin{definition}
	Let $\mathscr D$ be the Dynkin
	diagram associated to the root system $R$, and consider a subset
	$L\subset I$ of simple roots. We say that a fundamental co-weight $\omega_i^\vee$
	is \emph{essential relatively to $L$} if  each irreducible component
	of the graph $\mathscr 
	D\setminus\{\omega^\vee_i\}$ contains a root in $L$.
\end{definition}

% \begin{example}
% 	Consider the  Dynkin diagram
% 	$
% 	\dynkin[labels={1,2,3,4,5,6,7}] D{*o*o***}
% 	$,
% 	where the points marked as   $\dynkin A{o}$ correspond to the
% 	elements in $L=\{2,4\}$.  Then the
% 	fundamental co-weights essential relatively to $L$  are:
% 	$\{\omega^\vee_1,\omega^\vee_3,\omega^\vee_6,\omega^\vee_7\}$. 
% \end{example}

\begin{theorem}\label{th:DescriptionJlambda}
	If  $L\subset I$ then $J_L$ is the set of fundamental  co-weights that
	are essential relatively to $L$. 
\end{theorem}

\begin{proof}
	Under the  duality between $\Sigma_{R,L}$ and $\mathcal{WP}(\lambda)$
	(see Definition \ref{def:WeylPolytope} and  Proposition
        \ref{prop:WeylPolytope}) the elements of
        $\Prim{\sigma_{R,L}}$ 
	are the exterior normals of
         $C_\lambda(n-1)$. 
       % the facets of $\mathcal{WP}(\lambda)$  containing $\lambda$ (see Proposition
%	\ref{prop:WeylPolytope}).

        On the other hand,   by \cite[Theorem C]{kn:khare}
       an exterior normal of a facet in
           $C_\lambda(n-1)$ is in $\mathcal D$ if and only if the
           corresponding co-weight is essential (see \cite[Definition
           3.1]{kn:khare}). Since 
         $J_L$ 	corresponds to   the   facets   in $C_\lambda(n-1)$
        such that their exterior normal belongs to $\mathcal D$, the result follows.
\end{proof}

% In Proposition \ref{prop:faceSigmaLambda} above we describe the faces
% of the cone $\sigma_{R,L}$ in terms of the faces of $\mathcal D^\vee$;
We finish this section by presenting the notion of core of the cone
$\sigma_{R,L}$. This construction, that exhibits a relationship between the
 faces of $\mathcal D^\vee$  and the relative interior of
$\sigma_{R,L}$, will be useful in our characterization of the
Gorenstein Fano closures of generic orbits.

\begin{definition}\label{def:core}
	If $L\subset I$, we define the \emph{core} of $\sigma_{R,L}$, denoted
	by $\mathscr C(\sigma_{R,L})$, as the face of $\mathcal D^\vee$
	generated by the set $\{\omega^\vee_i\mathrel{:}i\in L^c\}$. 
\end{definition} 

The core of a cone $\sigma_{R,L}$ can easily be characterized by its
invariance properties: 

\begin{lemma}\label{lemma:core}
	If $L\subset I$,  then
	$\mathscr C(\sigma_{R,L})=(\mathcal
        D^\vee)^{W_L}=\bigcap_{w\in W_L}w \mathcal
        D^\vee=(\sigma_{R,L})^{W_L}.$ \qed
	\end{lemma}

% \begin{proof}
% 	The inclusions $\subset $ are obvious, as well as the first
% 	equality. If $v\in (\sigma_{R,L})^{W_L}$ 
% 	then, by definition, there exist $v'\in \mathcal D^\vee $ and $w\in
% 	W_L$ such that  $v=w\cdot v'$, and it  follows that
% 	$v=w^{-1}v=v'\in\mathcal D^\vee$.  \end{proof}

\begin{proposition}\label{prop:InteriorOfSigmaLambda}
	If  $L\subset I$ is such that $\sigma_{R,L}$ is 
	strictly convex, then 
	$\mathring{\mathscr C}(\sigma_{R,L})=({\mathring{\sigma}_{R,L}})^{W_L}$.	
\end{proposition}

\begin{proof} 	
  By
  % a classic result on the convex polyhedral cone and by using
  Proposition
  \ref{prop:dualcone},  $u\in \mathring{\sigma}_{R,L}$
  %belongs to the relative interior of $\sigma_{R,L}$
  if  and only if $\Dual{\beta}{u}>0$  for all $\beta\in (R^L)^+$.
		If  $v\in \mathring{\mathscr C}(\sigma_{R,L})$ then $
	v=\sum_{i\in L^c} a_i\omega^\vee_i$, with with $a_i>0$ for all $i\in
	L^c$. Thus, 
	$\Dual{\beta}{v}>0$ for all $\beta\in(R^L)^+$ 
          and therefore
 $\mathring{\mathscr C}(\sigma_{R,L})\subset \mathring{\sigma}_{R,L}$.
 We deduce from Lemma \ref{lemma:core} that
 $\mathring{\mathscr
   C}(\sigma_{R,L})\subset(\mathring{\sigma}_{R,L})^{W_L}$.

 Assume now that  $v\in
	(\mathring{\sigma}_{R,L})^{W_L}$. Then $v$ belongs to $({\mathcal
		D^\vee})^{W_L}\cap \mathring{\sigma}_{R,L}$ but,  by 
	Proposition \ref{prop:faceSigmaLambda}, $v$ does not belong to any of
	the facets of the cone ${(\mathcal D^\vee)}^{W_L}$. We apply
        again  Lemma \ref{lemma:core} and the result follows.  
\end{proof}

\subsection{A criteria for $X_{R,L}$ to be Gorenstein Fano}\ %

In this section  we characterize
when the  closure of a  generic orbit  is a Gorenstein
  Fano toric variety in 
terms  of the combinatorial properties of the core of the associated
cone  $\sigma_{R,L}$.

\begin{definition}\label{def:nlambda}
	We  define
	$\mathcal F_L$ as the convex hull of the set $\Prim{\sigma_{R,L}}$ and
	$\mathcal P_{R,L}$ as the convex hull of
        $\Prim{\Sigma_{R,L}}$.

        Clearly,    $\dim_\aff\mathcal P_{R,L}=n$ and 
	$\dim_\aff(\mathcal F_L)$ is $n$ or $n-1$ --- in this last case $\mathcal F_L$
	is a  facet of $\mathcal P_{R,L}$. When  $\mathcal F_L$ is a
        facet, we identify the exterior normal of $\mathcal F_L$ (see
        Definition \ref{defn:latpol}) with the unique
        element  $n_L\in (\Lambda_R^\vee)_{\QQ}$   such that
        $-\varphi_L(v)=(n_L,v)$.
  If  $\mathcal F_L$ is $n$-dimensional, we set $n_L=0$.
\end{definition}

\begin{remark}\label{rem:nlinv}\label{rem:nlk-1}
 (1) Notice that $n_L$ is $W_L$-invariant, since
	$\Prim{\sigma_{R,L}}=W_L\cdot J_L$.

        \noindent (2) In particular, if $L^c=\{\omega_i^\vee\}$, then 
        $n_L=a\omega_i^\vee$ --- this easy remark will simplify some
        calculations.

\noindent (3) Using the $W$-invariance of $n_L$, we can characterize
the affine dimension of  	$\Conv\bigl(\Prim{\sigma_{R,L}}\bigr)$  as follows:
	
	\begin{enumerate}[(i)]
\item 	If there exists a  $W_L$-invariant element   $n_L\in (\Lambda_R^\vee)_{\QQ}$ such that $ 
	\Scalar{n_L}{\nu}=1$ for all $\nu\in J_L$, then
	$\Conv\bigl(\Prim{\sigma_{R,L}}\bigr)$ is $(n-1)$-dimensional with
	exterior normal $n_L$.

	\item If such an element does not exist,  then
	$\Conv\bigl(\Prim{\sigma_{R,L}}\bigr)$ is $n$-dimensional.
\end{enumerate}

\end{remark}

\begin{theorem}\label{th:FinalCriteria}
	Let   $X_{R,L}$ be the closure of a generic orbit. Then $X_{R,L}$ is
	$\QQ$-Gorenstein Fano if and only if $n_L\in \mathring{\mathscr
		C}(\sigma_{R,L})$. 
Moreover, if $X_{R,L}$ is $\QQ$-Gorenstein Fano then 
	\[
	j_{_{X_{R,L}}}=\min\{j\in \NN^*\mathrel{ : }\forall v\in\Lambda_R^\vee,\,\Scalar{jn_L}{v}\in\ZZ\}.
	\]
%	where $j_{_{X_{R,L}}}$ is the Gorenstein index of $X_{R,L}$.
\end{theorem}

\begin{proof}
	Since $X_{R,L}=X_{-\Sigma_{R,L}}\cong X_{\Sigma_{R,L}}$ (see
	Definition \ref{def:FanRoot} and remark \ref{rem:prodfan}), it
	suffices to prove the assertion for $X_{\Sigma_{R,L}}$.
	
	Assume that
	$n_L\in \mathring{\mathscr 
		C}(\sigma_{R,L})$; then 
	  $\dim_\aff(\mathcal F_L)=(n-1)$ and, by
	Proposition \ref{prop:InteriorOfSigmaLambda},  $n_L\in
	(\mathring{\sigma}_{R,L})^{W_L}$. 
	If  $\sigma\in \Sigma_{R,L}(n)$ then   $\sigma=w\cdot
	\sigma_{R,L}$ for some $w\in W$, Therefore,  
	$
	\Prim{\sigma}=w\cdot \Prim{\sigma_{R,L}}= w W_L\cdot J_L$ and it follows that 
	\[
	w\cdot \mathring{\mathscr C}(\sigma_{R,L})=
	w\cdot    ({\mathring{\sigma}_{R,L}})^{W_L}\subset {w\cdot
		\mathring{\sigma}_{R,L}}= \mathring{\sigma}. 
	\]
	In particular,  $\operatorname{Conv}\bigl(\Prim{\sigma}\bigr)$ is
	$(n-1)$-dimensional, with exterior normal 
	$w\cdot n_L\in \mathring{\sigma}$.
	
	By Proposition \ref{prop:fanogeneral}, it remains to prove
	that $\Scalar{w\cdot n_L}{v}<1$ for all $v\in \Prim{\Sigma_{R,L}}\setminus
	\Prim{w\cdot \sigma_{R,L}}$. By the $W$-invariance
	of $\Sigma_{R,L}(n)$ and $\Prim{\Sigma}$ it suffices to prove that
	$\Scalar{ n_L}{v}<1$ for all $v\in \Prim{w\cdot \sigma_{R,L}}\setminus 
	\Prim{ \sigma_{R,L}}$, where $w\in W$ is such that $(w\cdot
	\sigma_{R,L})\cap\sigma_{R,L}$ is a common $(n-1)$-dimensional
	face.
	
	Moreover, by the $W_L$-invariance of $\sigma_{R,L}$, we can
	assume that $w=s_{\alpha_i^\vee}$, where $i\in L$ is such that
	$(\alpha_i^\vee)^\perp$ is the support hyperplane of a
	$(n-1)$-dimensional face of $\sigma_{R,L}$.
	In this case $v=s_{\alpha_i^\vee}\cdot \nu$ for some $\nu\in
	\Prim{\sigma_{R,L}}\setminus (\alpha_i^\vee)^\perp$ and we have that
	\[
	\Scalar{n_L}{s_{\alpha_i^\vee}\cdot  \nu}=  \Scalar{n_L}{ \nu} -
	\Scalar{n_L}{
		2\frac{\Scalar{\scriptstyle \alpha_i^\vee}{\scriptstyle \nu}}{\Scalar{\scriptstyle \alpha_i^\vee}{\scriptstyle \alpha_i^\vee}}\alpha_i^\vee}=1
	- \Scalar{n_L}{
		2\frac{\Scalar{\scriptstyle \alpha_i^\vee}{\scriptstyle \nu}}{\Scalar{\scriptstyle \alpha_i^\vee}{\scriptstyle \alpha_i^\vee}}\alpha_i^\vee}. 
	\]
	
	Since $\nu\in \sigma_{R,L}\setminus (\alpha_i^\vee)^\perp$ and  $n_L\in \mathring{\sigma_{R,L}}$, then 
	$\Scalar{\alpha_i^\vee}{\nu}>0$ and
	$\Scalar{n_L}{\alpha_i^\vee}>0$, and the assertion follows.

	Conversely,  assume now $n_L\notin \mathring{\mathscr
		C}(\sigma_{R,L})$; since  $n_L$ is
	$W_L$-invariant, it follows  that  $n_L=\sum_{i\in L^c}
	a_i\omega^\vee_i$.  If $n_L=0$ then
	$\operatorname{Conv}\bigl(\Prim{\sigma}\bigr)$ is $n$-dimensional and
	it follows from Proposition \ref{prop:fanogeneral} that $X_{R,L}$ is
	not $\QQ$-Gorenstein Fano. If $n_L\neq 0$,  there exists $i_0\in
	L^c$ such that $a_{i_0}\leq 0$. Since $i_0\in L^c$, it follows that
	$s_{\alpha_{i_0}^\vee}\notin W_L$ and therefore $s_{\alpha_{i_0}^\vee}\cdot \sigma_{R,L}
	\in \Sigma_{R,L}(n)\setminus \sigma_{R,L}$.  Let
	$\nu\in\Prim{\sigma_{R,L}}\setminus (\alpha_{i_0}^\vee)^\perp$. Then
	$\Scalar{n_L}{\nu}=1$ and $\Scalar{\nu}{\alpha_{i_0}^\vee}>0$. It follows
	that
{\small	\[
\scriptstyle	\Scalar{n_L}{s_{\alpha_{i_0}^\vee}\cdot  \nu}=  
	\Scalar{s_{\alpha_{i_0}^\vee}\cdot n_L}{ \nu} = 
	\Scalar{n_L}{\nu} -
	\Scalar{2\frac{\Scalar{\scriptstyle \alpha_{i_0}^\vee}{\scriptstyle n_L}}{\Scalar{\scriptstyle \alpha_{i_0}^\vee}{\scriptstyle \alpha_{i_0}^\vee}}\alpha_{i_0}^\vee}{\nu}=
		1- \Scalar{2a_{i_0}\frac{\Scalar{\scriptstyle \alpha_{i_0}^\vee}{\scriptstyle \omega_{i_0}^\vee}}{\Scalar{\scriptstyle \alpha_{i_0}^\vee}{\scriptstyle \alpha_{i_0}^\vee}}\alpha_{i_0}^\vee}{\nu}\geq1. 
	\]}

      Since $s_{\alpha_{i_0}^\vee}\cdot  \nu\in \Prim{s_{\alpha_{i_0}^\vee}\cdot \sigma_{R,L}}\setminus 
	\Prim{ \sigma_{R,L}}$, it follows from Proposition \ref{prop:fanogeneral} that $X_{\Sigma_{R,L}}$
	and hence  $X_{R,L}$ are not
	$\QQ$-Gorenstein Fano.

  The last assertion follows straightforward from Proposition
        \ref{prop:fanogeneral} --- indeed, recall that $W\cdot n_L$ is  set of
        exterior normals of the facets of 
  $\operatorname{Conv}\bigl(\Prim{\Sigma}\bigr)$.   	

\end{proof}

The following  (very)  easy  remark   will be used
several times in our classification of the Gorenstein Fano closures of
generic orbits.

\begin{lemma}\label{lemma:DiscardAlmostAll}
	Let $b_{i,j,k}=\Scalar{\omega^\vee_i}{\omega^\vee_{j}-\omega^\vee_{k}}$, $i,j,k\in I$. Assume that  $L\subset I$ is such  that there exist $j,k\in J_L$ with 
	$b_{i,j,k}\geq 0$ for all $i\in L^c$, and $\sum_{i\in L^c} b_{i,j,k}^2\neq
	0$. Then $n_L \notin \mathring{\mathscr
		C}(\sigma_{R,L}) $.
	
\end{lemma}

\begin{proof}
	If $n_L\in \mathring{\mathscr
		C}(\sigma_{R,L})$, then   $n_L=\sum_{i\in L^c}
	a_i\omega_i^\vee$, with $a_i>0$. 
	Since
	$\omega^\vee_{j}-\omega^\vee_{k}\in \bigl\langle
	\Prim{\sigma_{R,L}}\bigr\rangle_\aff$ (see Proposition
	\ref{prop:DescriptionNormal}), it follows that 
	$
	0< \sum_{i\in L^c}a_ib_{i,j,k}=\Scalar{n_L}{\omega^\vee_{j}-\omega^\vee_{k}}=0,
	$
	and we obtain a contradiction.
\end{proof}

\section{$\QQ$-Gorenstein Fano generic closures}
\label{section:table}

\begin{theorem}\label{th:Classification}
	Let $G$ be a simple affine algebraic group of root type $R$ and
	$L\subsetneq I$ a proper subset of the set of simple roots $I$;
	let $T$ be a maximal torus and $T\subset P$ be the parabolic subgroup
	associated to $L$. Table \ref{TabFinal}  on page
        \pageref{pagtable} gives a complete list
        of all the closures $X_{R,L}\subset 
	G/P$ of a generic $T$-orbit that are $\QQ$-Gorenstein Fano,
	Gorenstein Fano and smooth Fano varieties --- e.g. if
          the table indicates that $X_{R,L}$ is $\QQ$-Gorenstein Fano,
          then $X_{R,L}$ is  $\QQ$-Gorenstein Fano but not Gorenstein Fano.

        In the third column,
	we draw the  Dynkin diagram of $R^\vee$; the subscripts indicate the
	number of the corresponding simple roots in $I$, the elements of $L$
	are the roots drawn in black, the superscript $J_L$ over a root $\alpha_i$
	indicates that $\omega_i^\vee\in J_L$, the set of essential
	co-weights
        %--- recall that $\Prim{\sigma_{R,L}}=W_L\cdot J_L$ and that
	% $\Prim{\Sigma_{R,L}}=W\cdot J_L$.
        In the fourth column we indicate
	the corresponding geometry ---  the Gorenstein index of the
	$\QQ$-Gorenstein Fano variety $X_{R,L}$ is denoted by  $j$. Finally, in the fifth column we exhibit
	the exterior normal  $-\varphi_L\in (\Lambda_R)_\QQ$   ---
        recall that 
	$\langle -\varphi_L, v\rangle=\Scalar{n_L}{v}$ for all $v\in (\Lambda_R^\vee)_\QQ$.
\end{theorem}

\begin{remark} (1) The description given in Table \ref{TabFinal} is
	modulo automorphisms of the root system; for 
	example, the varieties $X_{A_n,I\setminus\{1\}}$ and
	$X_{A_n,I\setminus\{n\}}$ are isomorphic (and both Fano) and
        the table
	only exhibits $X_{A_n,I\setminus\{1\}}$. In the same
	spirit, the conditions given on the rank are established in order  to
	avoid repetition.      
	
	\noindent (2) In the proof of Theorem \ref{th:Classification} we
	exhibit $n_L$ when $X_{R,L}$ is $\QQ$-Gorenstein Fano; in order to
	compute $\varphi_{L}$  we use the matrix
	$\big(\Dual{\omega_i}{\omega^\vee_j}\big)_{i,j\in I}$, 
	which is the inverse of the Cartan Matrix (see \cite{kn:OniVinberg}
	for explicit calculations).
\end{remark}

% \newpage

{
\begin{figure}
   \captionof{table}{\label{TabFinal}
  }\label{pagtable}
	\newcommand{\dsp}{\displaystyle}
	\renewcommand\arraystretch{2.5}
	{\footnotesize
		\setlength{\doublerulesep}{0pt}
		\begin{tabular}{||c|l|c|l|c||}
			\hline
			\hline
			$\operatorname{type}$ &rank& Dynkin($R^\vee$), $L$ and $J_L$&
			Geometry&$-\varphi_L\in(\Lambda_R)_\QQ$ \\
			\hline 
			\hline 
			%%%%%%%%%%%%%%%%%%%%%%%%%%%%%%%%%%%%%%%%%%%%%%%%
			\multirow{4}{*}{$A_n$}& $n\geq
                                                1$&$\dynkin[labels*={,,J_L},labels={1,2,n}]A{o*.*}
                                                    $& Smooth,
                                                       Fano&$
                                                             (n+1)\omega_1$\\
			\cline{2-5}
			&$n\geq 2$    &$\dynkin[labels*={J_L,J_L,J_L,J_L},labels={1,2,n-1,n}]A{o*.*o}$& Gorenstein Fano&$\omega_1+\omega_n$\\
			\cline{2-5}
			&  $n\geq 3$,  odd &$\dynkin[labels*={J_L,,,,J_L},labels={1,\ ,\frac{n-1}{2},\ ,n-1}]A{*.*o*.*}$&Gorenstein Fano&$2\omega_{\frac{n+1}{2}}$\\
			\cline{2-5}
			& $n\geq 4$,  even &$\dynkin[labels*={J_L,,,,,J_L},labels={1,\ ,\frac n2,\frac n2+1,\ ,n}]A{*.*oo*.*}$&	Smooth, Fano&  $\scriptstyle(n+1)(\omega_{_{\frac n2}}+\omega_{_{\frac n2 +1}})$\\
			\hline
			\hline
			%%%%%%%%%%%%%%%%%%%%%%%%%%%%%%%%%%%%%%%%%%%%%%%%%%%%%
			\multirow{2}{*}{$B_n$}&  \multirow{3}{*}{$n\geq 2$}&$\dynkin[labels*={,,,J_L},labels={1,2,n-1,n}]C{o*.**}$&Gorenstein Fano&$\omega_1$\\
			\cline{3-5}
			&  &$\dynkin[labels*={J_L,,,},labels={1,2,n-1,n}]C{**.*o}$&Smooth, Fano&$2\omega_n$\\
			\hline
			\hline
			%%%%%%%%%%%%%%%%%%%%%%%%%%%%%%%%%%%%%%%%%%%%%%%%%%%
			\multirow{3}{*}{$C_n$}& \multirow{3}{*}{$n\geq 3$}&$\dynkin[labels*={,,,J_L},labels={1,2,n-1,n}]B{o*.**}$&Gorenstein Fano&$\omega_1$\\
			\cline{3-5}
			&&$\dynkin[labels*={J_L,,,,J_L},labels={1,2,\ ,n-1,n}]B{*o*.**}$&	 Gorenstein Fano&$\omega_2$\\
			\cline{3-5}
			&&$\dynkin[labels*={J_L,,,},labels={1,2,n-1,n}]B{**.*o}$&
			$n$
			even:
			Gor. Fano &$\omega_n$\\
			&&&$n$ odd: $\QQ$-G.F., $\scriptstyle j=2$&\\
			
			\hline
			\hline
			
			%%%%%%%%%%%%%%%%%%%%%%%%%%%%%%%%%%%%%%%%%%%%%%%%%%%%%%%%%%%
			\multirow{2}{*}{$D_n$}&\multirow{2}{*}{ $n\geq 4$}&$\dynkin[labels*={,,,J_L,J_L},labels={1,2,,n-1,n}]D{o*.***}$&Gorenstein Fano&$2\omega_1$\\
			\cline{3-5}
			&	&$\dynkin[labels*={J_L,,,J_L,J_L},labels={1,2,,n-1,n}]D{*o.***}$&Gorenstein Fano&$\omega_2$\\
			\hline
			\hline
			
			%%%%%%%%%%%%%%%%%%%%%%%%%%%%%%%%%%%%%%%%%%%%%%%%%%%%%%%%
			$E_6$&6&$\dynkin[labels*={J_L,,,,,J_L},labels={1,2,3,4,5,6}]E{*o****}$& Gorenstein Fano&$\omega_2$\\	
			\hline
			\hline
			%%%%%%%%%%%%%%%%%%%%%%%%%%%%%%%%%%%%%%%%%%%%%%%
			\multirow{2}{*}{$F_4$}&\multirow{2}{*}{4}&$\dynkin[labels*={,,,J_L},labels={1,2,3,4},reverse
                                                                   arrows]F{o***}$&
                                                                                    $\QQ$-Gor. Fano, $j=2$&$\frac{\omega_1}{2}$\\
			\cline{3-5}
			&&$\dynkin[labels={1,2,3,4},labels*={J_L,,,},reverse arrows]F{***o}$&Gorenstein Fano&$\omega_4$\\
			\hline
			\hline
			
			%%%%%%%%%%%%%%%%%%%%%%%%%%
			\multirow{2}{*}{$G_2$}&\multirow{2}{*}{2}&$\dynkin[labels={1,2},labels*={J_L,},reverse arrows]G{*o}$&Smooth, Fano&$\omega_2$\\
			\cline{3-5}
			&&$\dynkin[labels={1,2},labels*={,J_L},reverse
                           arrows]G{o*}$&$\mathbb Q$-Gor. Fano, $j=3$&$\frac{\omega_1}{3}$\\
			\hline
			\hline
			
		\end{tabular}
	}
\end{figure}}

%\newpage 

\section{Proof of theorem \ref{th:Classification} }
\label{sec:explicalc}

\subsection{Explicit calculations for  ranks $1$ and $2$}\ %

The  classification for ranks $1$ and $2$ is  done by direct examination.

\medskip

\noindent {\bf Explicit calculations for $\mathbf G$ of type $\mathbf{A_1}$.}

  Clearly
$\sigma_{A_1,\{1\}}= \QQ^+$ and therefore
$X_{A_1,\emptyset}=\operatorname{SL}_2(\Bbbk)/B\cong \PP^1$ is a smooth Fano
variety.
%--- recall that  $\sigma_{R,I}=(\Lambda_R)^\vee_\QQ$ for all
%root system $R$, see Proposition \ref{prop:sigma.convex}.

\medskip

% \noindent {\bf Explicit calculations for $\mathbf G$ of type $\mathbf{A_1\times A_1}$.}

%  By Proposition \ref{prop:sigma.convex}, the only generic closure is $X_{A_1\times A_1,\emptyset}=G/B\cong
% \PP^1\times\PP^1$; this is  clearly
% $X_{A_1\times A_1,\emptyset}$ a smooth Fano variety. See also Remark
% \ref{rem:prodfan}.   

% It follows from Remark \ref{rem:prodfan} that  that is
%a smooth Fano variety
%--- this case is covered by  Remark \ref{rem:prodfan},
%but we include it for the sake of completeness ``in small ranks''.

% In this case the Dynkin diagram is:
% \[
% \dynkin[labels={1}] A1\ \ \ \dynkin[labels={2}] A1
% \]

%  To have the hypothesis of the proposition \ref{prop:sigma.convex} verified, we have to choose $L=I=\{1,2\}$; then $X_{R,L}$ is isomorphic to $\PP^1\times \PP^1$ and is a Fano variety.  

% \begin{figure}[h!]
% \centering
% \input{A1A1.pdf_t}
% \caption{}
% \label{picture:A1A1}
% \end{figure}

\medskip

\noindent {\bf Explicit calculations for $\mathbf G$ of type $\mathbf{A_2}$.} \ %

\begin{figure}[h!]
	\centering
	\subfloat[$L=\{2\}$\label{picture:A2_1}]{\resizebox{.28\textwidth}{!}{
			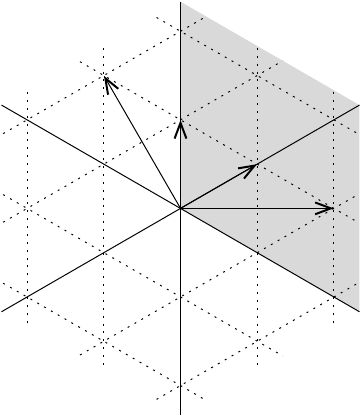}}
	\hspace*{0.2\textwidth}
	\subfloat[$L=\emptyset$ \label{picture:A2_12} ]{\resizebox{.28\textwidth}{!}{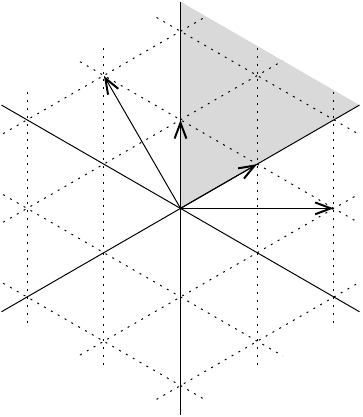}}
	\caption{$G$ is of type $A_2$;  the
cone $\sigma_{R,L}$ is drawn in   gray.  
}
	\label{picture:A2}
\end{figure}

As follows from  figures
\ref{picture:A2_1} and \ref{picture:A2_12}, $\Sigma_{A_2,L}$ is
a smooth Fano complete fan  for 
all $L\subsetneq I =\{1,2\}$ --- recall that $R\cong R^\vee$. The
varieties  $X_{A_2,\{1\}}$ and 
$X_{A_2,\{2\}}$ are isomorphic to $\PP^2$ with the canonical action of
$T=(\Bbbk^*)^2\subset \PP^2$,  and $X_{A_2,\emptyset}$ is
isomorphic to the blowing up of  three generic points in $\PP^2$.

\medskip

\noindent {\bf Explicit calculations for $\mathbf G$ of type
	$\mathbf{B_2}$ and $\mathbf{C_2}$.}% \label{explicit:B2} \ %

  Since   $C_2=(B_2)^\vee$,  it suffices to describe the case
    $R=B_2$, which has   Dynkin diagram $\dynkin[edge length=0.8cm,label,arrow width=1.5mm,] B2$.

% {\color{blue}    Since   $C_2=(B_2)^\vee$,  it suffices to describe the case
%     $R=B_2$, which has  recall that the Dynkin diagram of $B_2$ is $
% \dynkin[edge length=0.8cm,label,arrow width=1.5mm,
% ] B2$. }

% Since $C_2$ is the dual root system of $B_2$,  it suffices to describe the case $B_2$.
% the Dynkin diagram of $R=B_2$ is $
% \dynkin[edge length=1cm,label,arrow width=1.5mm,
% ] B2$ whereas the  Dynkin diagram of $R^\vee=C_2$ is $\dynkin[edge length=1cm,label,arrow width=1.5mm,
% ]  C2$. Therefore, it suffices to describe the case $B_2$.

\begin{figure}[h!]
	\centering
	\subfloat[$L=\{1\}$\label{picture:B2_2}]{\resizebox{.25\textwidth}{!}{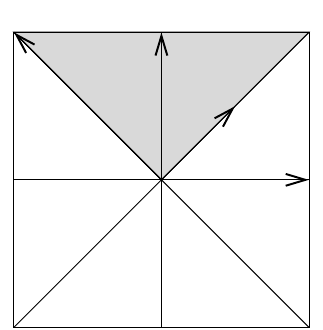}}
	\hspace*{0.25\textwidth}\subfloat[$L=\{2\}$\label{picture:B2_1}]{\resizebox{.25\textwidth}{!}{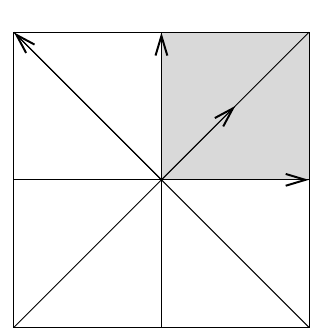}}	
	\caption{$G$ is of type  $B_2$. The cone
          $\sigma_{R,L}$ is drawn in gray.  
}
	\label{picture:B2}	
\end{figure}

If $L=\emptyset$, then $n_{\emptyset}=\omega_1^\vee$ and therefore
$X_{B_2,\emptyset}$ is not $\QQ$-Gorenstein Fano (see Figure
\ref{picture:B2_2}). If follows that $X_{C_2,\emptyset}$ is not
$\QQ$-Gorenstein Fano --- in this case $n_\emptyset=\omega_2^\vee$.

If $L=\{1\}$ then $n_{\{1\}}=\omega^\vee_2$ and  $X_{B_2,\{1\}}$ is a
smooth Fano variety ---  $X_{B_2,\{1\}}$ is isomorphic to $\PP^1\times\PP^1$,
see figure \ref{picture:B2_2}. It follows
that $X_{C_2,\{2\}}\cong \PP^1\times\PP^1$. 

If  $L=\{2\}$, then $n_{B_2,\{2\}}=\omega_1^\vee$  and $X_{B_2,
	\{2\}}$ and $X_{C_2,\{1\}}$ are Gorenstein Fano varieties  (see figure
\ref{picture:B2_1}) --- notice that
these varieties are  not smooth.

\medskip

\noindent {\bf Explicit calculations for $\mathbf G$ of type $\mathbf{G_2}$.}\ %
\label{exam:G2}

Recall that $G_2$ has   Dynkin diagram  $ \dynkin[edge
length=0.8cm,label,arrow width=1.5mm]G2$, with $G_2\cong G_2^\vee$.

% Recall that  the matrix of the $W$-invariant scalar product on
% $(\Lambda_R^\vee)_\QQ$ is $\left( \begin{smallmatrix}
% 2&3\\
% 3&6\\
% \end{smallmatrix} \right)$.

\begin{figure}[h!]
	\centering
	\subfloat[$L=\{1\}$\label{picture:G2_2}]{\resizebox{.28\textwidth}{!}{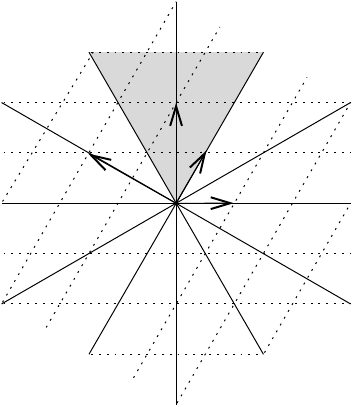}}
	\hspace*{.18\textwidth}
	\subfloat[$L=\{2\}$\label{picture:G2_1}]{\resizebox{.28\textwidth}{!}{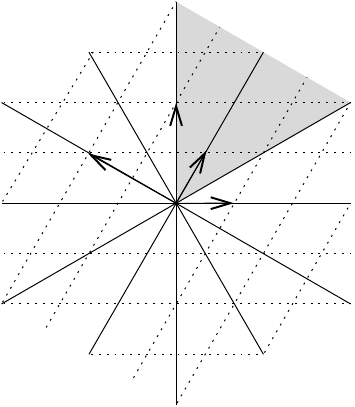}}	
	\caption{ $G$ is of type  $G_2$}
	\label{picture:G2}
\end{figure}

Since $\Conv\bigl(\Prim{\sigma_{G_2,\emptyset}}\bigr)$ is not a proper
face of $\Conv\bigl(\Prim{\Sigma_{G_2,\emptyset}}\bigr)$ (see figure
\ref{picture:G2_2}), it follows that $X_{G_2,\emptyset}$ is not $\QQ$-Gorenstein Fano.

An inspection of   figure \ref{picture:G2_2} shows that
$X_{G_2,\{1\}}$ is a smooth Fano variety --- in fact, $X_{G_2,\{1\}}\cong
X_{A_2,\emptyset}$, the blowing up of three points in 
$\PP^2$. 
Also by  inspection of figure  \ref{picture:G2_1}, we have that
$X_{G_2,\{2\}}$ is  $\QQ$-Gorenstein Fano of index $3$. 

% \begin{example}\label{ex:exoticfan}
% Let $G$ be of type $G_2$, and consider $\Sigma$ as in the picture
% \ref{picture:G2_excp}; clearly, $\Sigma$ is a smooth complete fan.  However,
% $\Sigma$ does not correspond to a generic $T$-orbit of an homogeneous
% space $G/P$. Thus, there exist Fano toric varieties whose maximal
% cones are union of Weyl chambers, but are not of the form
% $\Sigma_{R,L}$ for some subset $L$. 
% \begin{figure}[h!]
% 	\centering
% 	\input{G2_excp.pdf_t}
% 	\caption{}
% 	\label{picture:G2_excp}
% \end{figure}

% \end{example}

\subsection{Proof of theorem \ref{th:Classification} for rank $n\geq
  3$}\ %
\label{sec:explicalc2}

% Here we give the calculation for all type for rank $n\geq 3$.
% The  strategy is similar for all cases.
%  For each root system $R$,
% we describe for the dual root system $R^\vee$ a matrix of a  $W$-invariant scalar
% product in the basis of the fundamental co-weights. Notice
% that this matrix is a multiple of the symmetrized of the inverse of
% the Cartan Matrix associated to $R^\vee$ --- see for example
% \cite[p. 295]{kn:OniVinberg}. By Theorem \ref{th:FinalCriteria}, we
% need to determine when  $n_L$ belongs  to $ \mathring{\mathscr
% 	C}(\sigma_{R,L})$. Hence,  for each
% $L$, we use Proposition \ref{prop:faceSigmaLambda} in order to compute
% $J_L$ and then  we deduce $n_L$ (see Notation \ref{nota:JL} and
% Definition \ref{def:nlambda}) --- sometimes, we will assume that
% $X_{R,L}$ is $\QQ$-Gorenstein Fano and obtain a contradiction. Notice that for simplicity of
% computing, we reason on $L^c$ rather  than $L$.
 We follow a similar strategy for all cases: given $R$,
 (a multiple of) the transposed inverse of  
symmetrized  Cartan Matrix associated to $R^\vee$ determines a $W$-invariant
scalar product on $R^\vee$ (see for example
\cite[p. 295]{kn:OniVinberg}).  
 By Theorem \ref{th:FinalCriteria}, we
need to determine when  $n_L$ belongs  to $ \mathring{\mathscr
	C}(\sigma_{R,L})$ (see Definition \ref{def:nlambda}). Hence,  for each
$L$, we use Proposition \ref{prop:faceSigmaLambda} in order to compute
$J_L$ and then  we deduce $n_L$.

 % --- sometimes, we will assume that
 % $X_{R,L}$ is $\QQ$-Gorenstein Fano and obtain a contradiction.

 \begin{notation}
Notice that  we reason on $L^c$ rather  than on $L$. If  $L\subset
I=\{1,2,\ldots,n\}$, we denote by $m$ (resp. $M$) the 
minimum (resp. the maximum) of the set $L^c$.   
\end{notation}

\subsection{Explicit calculations for $G$ of type $A_n$, $n\geq 3$}\ %
%Here $R$ and $R^\vee$ are isomorphic, with   following Dynkin diagram: 
%\[
%\dynkin[edge length=1cm,labels={1,2,n-1,n}] A{}.
%\]

We choose the $W$-invariant scalar product  given by the  matrix:
% \[
% \begin{pmatrix}
% n&n-1&n-2&\ldots&\ldots&\ldots&2&1\\
% n-1&2.(n-1)&2(n-2)&\ldots&\ldots&\ldots&2.2&2\\
% \vdots&\vdots& &\ldots&\ldots&\ldots&\vdots&\vdots\\
% n-i+1&2.(n-i+1)&\ldots&i.(n-i+1)&i.(n-i)&\ldots&i.2&i\\
% \vdots&\vdots& &\ldots&\ldots&\ldots&\vdots&\vdots\\
% 1&2&3&\ldots&\ldots&\ldots&(n-1)2&n
% \end{pmatrix}
% \]
% \begin{footnotesize}
	\begin{equation}\label{eqn:scalarAn}
	 \left(\vcenter{\xymatrixcolsep{0.0001pc}\xymatrixrowsep{0.0001pc}\xymatrix{
			\scriptstyle n&\scriptstyle  n-1&\scriptstyle  n-2\ar@{.}[rrrr]&&&& \scriptstyle 2&\scriptstyle  1\\
		\scriptstyle 	n-1\ar@{.}[dd]&\scriptstyle  2(n-1)\ar@{.}[dd]&\scriptstyle  2(n-2)\ar@{.}[rrrr]&&&&\scriptstyle  2\cdot 2\ar@{.}[dd]&\scriptstyle  2\ar@{.}[dd]\\
			&&&&&&\\
\scriptstyle 			 n-i+1\ar@{.}[dd]& \scriptstyle 2(n-i+1)\ar@{.}[dd]\ar@{.}[rr]&&
\scriptstyle 			i(n-i+1)&\scriptstyle  i(n-i)\ar@{.}[rr]&&\ar@{.}[dd]
		\scriptstyle 	i\cdot
			2&\ar@{.}[dd]\scriptstyle  i\\
			&&&&&&\\
		\scriptstyle 	2& \scriptstyle 2\cdot 2\ \ar@{.}[rrrrr]&&&&&\scriptstyle  (n-1)2& \scriptstyle n-1\\
			\scriptstyle 1&\scriptstyle  2\ \ar@{.}[rrrrr]&&&&&\scriptstyle  (n-1)&\scriptstyle  n
	}}\right)
	\end{equation}
%\end{footnotesize}

%   \[
% \left(\vcenter{\xymatrixcolsep{0.001pc}\xymatrixrowsep{0.001pc}\xymatrix{
%     \scriptstyle n&\scriptstyle n-1&\scriptstyle n-2\ar@{.}[rrrr]&&&&\scriptstyle 2&\scriptstyle 1\\
%    \scriptstyle  n-1\ar@{.}[dd]&\scriptstyle 2(n-1)\ar@{.}[dd]&\scriptstyle 2(n-2)\ar@{.}[rrrr]&&&&\scriptstyle 2\cdot 2\ar@{.}[dd]&\scriptstyle 2\ar@{.}[dd]\\
%     &&&&&&\\
%    \scriptstyle
%    n-i+1\ar@{.}[dd]&\scriptstyle 2(n-i+1)\ar@{.}[dd]\ar@{.}[rr]&&\scriptstyle
%    i(n-i+1)&\scriptstyle i(n-i)\ar@{.}[rr]&&\ar@{.}[dd]\scriptstyle
%    i\cdot
%     2&\ar@{.}[dd]\scriptstyle i\\
%      &&&&&&\\
% \scriptstyle 1&\scriptstyle 2&\scriptstyle 3\ \ar@{.}[rrrr]&&&&\scriptstyle (n-1)2&\scriptstyle n
% }}\right)
% \]

\subsubsection{Cases $L^c=\{1\}$ and $L^c=\{n\}$}\ %

If $L^c=\{1\}$ then  $J_L=\{\omega^\vee_n\}$, and it follows from
Remark \ref{rem:nlinv} that  $n_L=\omega_1^\vee$. Hence,  $X_{A_n,L}$ is Gorenstein Fano (see Theorem \ref{th:FinalCriteria}). Moreover,  a direct
computation (see also \cite[Remark 4]{kn:voskly}) shows that
\[
\Prim{\sigma_{R,L}}=W_{L}\cdot (\omega^\vee_n)=\{\omega^\vee_n,\omega^\vee_{n-1}-\omega^\vee_{n},\omega^\vee_{n-2}-\omega^\vee_{n-1},\ldots,\omega^\vee_1-\omega^\vee_2\},
\]
and therefore  $X_{A_n,L}$ is a  smooth  Fano variety.

By symmetry of the Dynkin diagram, we deduce that $X_{A_n, I\setminus
	\{n\}}$ is also a smooth Fano variety.

\subsubsection{Case $L^c=\{i\}$, $ i\neq 1,n$} \ %

In this case $J_L=\{\omega^\vee_1,\omega^\vee_n\}$, and
$n_L=a\omega_i^\vee$ (see Remark \ref{rem:nlk-1}). If $n\neq 2i-1$
then
$\Scalar{\omega^\vee_1}{\omega^\vee_i}=\Scalar{\omega^\vee_n}{\omega^\vee_i}$,
and therefore $n_L=0$. It follows that the polytope $\mathcal F_L$
(see Definition \ref{def:nlambda}) is
$n$-dimensional and hence  $X_{A_n,L}$
is not $\QQ$-Gorenstein Fano (see Proposition \ref{prop:fanogeneral}).

If $n=2i-1$, then $n_L=\frac{\omega^\vee_i}{i}$ and
$\Scalar{n_L}{\omega^\vee_j}\in\ZZ$ for all $j\in I$ (see Matrix \eqref{eqn:scalarAn}). It follows that
$X_{A_n,I\setminus \{i\}}$ is Gorenstein Fano.
% --- indeed, notice that
% in this case  $n-i+1=i$ and the  scalar products are given by the
% $i^{\text{th}}$ column of the  matrix .
Finally, since
{\small 
\[
\begin{split}
 \Prim{\sigma_{R,L}}= & \ W_L\cdot \omega^\vee_1\cup W_L\cdot
\omega^\vee_n =\\
&\  \{\omega^\vee_1,\omega^\vee_2-\omega^\vee_1,\ldots,\omega^\vee_{i}-\omega^\vee_{i-1}\}\cup
\{\omega^\vee_n,\omega^\vee_{n-1}-\omega^\vee_{n},\ldots,
\omega^\vee_{i}-\omega^\vee_{i+1}\} 
\end{split} 
\]}
 is not a  simplicial set,  it
follows that $X_{A_n,L}$ is not smooth.

\subsubsection{Case $M=m+1$} \ %

In this case $J_L=\{1,n\}$. If $X_{A_n,L}$ is
$\QQ$-Gorenstein Fano, it follows from Theorem \ref{th:FinalCriteria} that $n_L=a\omega^\vee_m+b\omega^\vee_{m+1}$,
with $a,b\in\QQ_{>0}$. Hence, 
  \[
0= \Scalar{n_L}{\omega^\vee_1-\omega^\vee_n} = (a+b)(n-2m)+(a-b).
\]
% where the last equality follows from the fact that
% $\Scalar{\omega^\vee_{i}}{\omega^\vee_1-\omega^\vee_n}={n-2i+1}$.
% If follows that

We deduce that $n=2m$ and
$n_L=(\omega^\vee_m+\omega^\vee_{m+1})$. Moreover, 
{\small
\[
\begin{split}
& \Prim{\sigma_{R,L}}= \ W_L\cdot \omega^\vee_1\cup W_L\cdot
\omega^\vee_n =\\
&\quad  \quad \{\omega^\vee_1,\omega^\vee_2-\omega^\vee_1,\ldots,\omega^\vee_{m}-\omega^\vee_{m-1}\}\cup
\{\omega^\vee_n,\omega^\vee_{n-1}-\omega^\vee_{n},\ldots,
\omega^\vee_{m+1}-\omega^\vee_{m+2}\},  
\end{split} 
\]
}
\noindent and therefore $\Prim{\sigma_{R,L}}$ is  a  basis of the co-weight lattice. It
follows that the variety $X_{A_{2m},\{m,m+1\}^c}$ is smooth Fano.

\subsubsection{Case $m=1$ and $M=n$}\ %

In this case,  $J_L=I$. If $X_{A_n,L}$ is $\QQ$-Gorenstein Fano, then $n_L=\sum_{i\in
	L^c}a_i\omega_i^\vee$, $a_i>  0$, is defined up to a scalar by the equations
$\Scalar{n_L}{\omega^\vee_j-\omega^\vee_{j'}}=0$. Since  
$\Scalar{\omega^\vee_1+\omega^\vee_n}{\omega^\vee_j-\omega^\vee_{j'}}=0$ 
we deduce that  $n_L$ is proportional to
$\omega^\vee_1+\omega^\vee_n$, and therefore:

\noindent (i) If $L^c\neq \{1,n\}$, then $X_{A_n,L}$ is not $\QQ$-Gorenstein Fano. 

\noindent  (ii) If $L^c=\{1,n\}$, then
$n_L=\frac{\omega^\vee_1+\omega^\vee_n}{n+1}$ and  therefore 
$X_{A_n,\{2,\dots,n-1\}}$ is Gorenstein Fano. Since $
\Prim{\sigma_{A_n,\{2,\dots,n-1\}}}= W_L\cdot\{\omega^\vee_j\ |\ j\in I\}$
is clearly not simplicial, we deduce that $X_{A_n, \{1,n\}^c}$ is not
smooth Fano.

\subsubsection{Remaining cases}\ %

It remains to study the cases  (i)  $\#L^c\geq 2$, with $m\neq M-1$
and $M<n$, and (ii) $\#L^c\geq 2$, with  $m\neq
M-1$ and $ 1<m$.

By symmetry, it suffices to study the case  (i); in this case  $J_L=\{1,m+1,\ldots, M-1,n\}$. 
Assume that $n_L=\sum_{i\in L^c}a_i\omega^\vee_i$,  $a_i> 0$; since 
$n$ and $M-1$ belong to $J_L$, we deduce  that 
\[\sum_{i\in L^c}a_i\Scalar{\omega^\vee_i}{\omega^\vee_{M-1}-\omega^\vee_{n}}=0.
\]

The coefficients of the $i$-th row of Matrix \eqref{eqn:scalarAn}  being a strict unimodal sequence with peak at the $i$-place, we deduce that  
$\Scalar{\omega^\vee_i}{\omega^\vee_{M-1}-\omega^\vee_{n}}>0$ for all
$i$ such that $i\leq M-1$. Since  $M<n$, we have that
$\Scalar{\omega^\vee_{M}}{\omega^\vee_{M-1}-\omega^\vee_{n}}>0$. Hence,
since $L\subset \{m,\ldots,M\}$, we obtain a contradiction (using   Lemma  
\ref{lemma:DiscardAlmostAll}); it follows  that  $X_{A_n,L}$ is not $\QQ$-Gorenstein Fano.

\subsection{Explicit calculations for $G$ of type $B_n$, $n\geq 3$}\ %
% Recall that  $B_n$ has Dynkin diagram 
% $
% \dynkin[labels={1,2,n-1,n}] B{**.**}
% $
% and    $(B_n)^\vee$  is of type $C_n$, with Dynkin diagram
% $
% \dynkin[labels={1,2,n-1,n}] C{**.**}$.

% Recall that a $W$-invariant scalar product is
% given by the following matrix in the
% basis of the fundamental co-weights (relatively to $C_n=(B_n)^\vee$):

Recall that  $B_n$ has Dynkin diagram 
$
\dynkin[labels={1,2,n-1,n}] B{**.**}
$
and    $(B_n)^\vee=C_n$.
We choose the  $W$-invariant scalar product on $(\Lambda_{B_n}^\vee)_{_\QQ}=(\Lambda_{C_n})_{_\QQ}$
given by the matrix
%\begin{footnotesize}
	\begin{equation}\label{eqn:scalarCn}
	\left(\vcenter{\xymatrixcolsep{0.0001pc}\xymatrixrowsep{0.0001pc}\xymatrix{
		\scriptstyle 	1&\scriptstyle 1\scriptstyle &\scriptstyle 1\ar@{.}[rrr]&&&\scriptstyle 1\ar@{.}[rr]&&\scriptstyle 1&\scriptstyle 1\\
			\scriptstyle 1\ar@{.}[dd]&\scriptstyle 2\ar@{.}[dd]&\scriptstyle 2\ar@{.}[rrr]\ar@{.}[dd]&&&\scriptstyle 2\ar@{.}[dd]\ar@{.}[rr]&&\scriptstyle 2\ar@{.}[dd]&\scriptstyle 2\ar@{.}[dd]\\
			&&&&&&&&\\
			\scriptstyle 1\ar@{.}[dd]&\scriptstyle 2\ar@{.}[dd]&\scriptstyle 3\ar@{.}[rrr]\ar@{.}[dd]&&&\scriptstyle i\ar@{.}[dd]\ar@{.}[rr]&&\scriptstyle i\ar@{.}[dd]&\scriptstyle i\ar@{.}[dd]\\
			&&&&&&&\\
			\scriptstyle 1&\scriptstyle 2&\scriptstyle 3\ar@{.}[rrr]&&&\scriptstyle i\ar@{.}[rr]&&\scriptstyle n-1&\scriptstyle n
	}}\right)
	\end{equation}
%\end{footnotesize}

\subsubsection{Case $\#L^c=1$}\ %
\label{case:Bn-1}

\noindent (i) If  $L^c=\{i\}$ with $i\neq 1,n$,  then
$J_L=\{1,n\}$ and $n_L=a\omega^\vee_i$. Since
$\Scalar{\omega_i^\vee}{\omega^\vee_1-\omega^\vee_n}=1-i$, it
follows that  $n_L=0$ and  $X_{B_n,\{i\}^c}$ is not $\QQ$-Gorenstein Fano.  \medskip

\noindent (ii) If $L^c=\{1\}$, then $J_L=\{n\}$. Since
$\Scalar{\omega^\vee_{1}}{\omega^\vee_{n}}=1$, it follows  that
$n_L=\omega_1^\vee$ and therefore 
$X_{B_n,\{2,\dots,n\}}$ is Gorenstein Fano. Since   $\#
W_L\cdot\omega^\vee_n>n$, then   $\sigma_{B_n,\{2,\dots,n\}}$
is not simplicial and we deduce that  $X_{B_n,\{2,\dots,n\}}$ is not
smooth.\medskip 
% If $L^c=\{1\}$, then $J_L=\{n\}$, and 
% $X_{B_n,\{2,\dots,n\}}$ is Gorenstein Fano --- since
% $\Scalar{\omega^\vee_{1}}{\omega^\vee_{n}}=1$, see  Theorem
% \ref{th:FinalCriteria}. It is easy to see that  $\#
% W_L\cdot\omega^\vee_n>n$, therefore,   $\sigma_{B_n,\{2,\dots,n\}}$
% is not simplicial and the variety $X_{B_n,\{2,\dots,n\}}$ is not
% Fano.\\ 

\noindent (iii) If  $L^c=\{n\}$, then $J_L=\{1\}$ and it  follows
that $n_L=\omega^\vee_n$. Hence, 
$X_{B_n,\{1,\dots,n-1\}}$ is  Gorenstein Fano. Moreover, it is easy to
see that $W_L\cdot J_L=
W_L\cdot\omega_{1}^\vee=\{\omega_1^\vee,\omega_2^\vee-\omega_1^\vee,\ldots,\omega_{n}^\vee-\omega_{n-1}^\vee\}$
(see for example  \cite[Remark (4)]{kn:voskly}) and therefore the
variety $X_{B_n,\{1,\dots, n-1\}}$ is smooth Fano.

\bigskip

\subsubsection{Case $\#L^c>1$}\label{sec:BnL>2}\ %

It is clear that    $\{1,n\}\subset J_L$. Since
$\Scalar{\omega^\vee_{i}}{\omega^\vee_n-\omega^\vee_{1}} = i-1$ for
all $i\in
I$, we deduce from     Lemma
\ref{lemma:DiscardAlmostAll} that $X_{B_n,L}$ is not $\QQ$-Gorenstein Fano.

\subsection{Explicit calculations for $G$ of type $C_n$}\ %

Recall that $C_n$ has   Dynkin diagram
$
\dynkin[labels={1,2,n-1,n}] C{**.**}$. We choose the  $W$-invariant
scalar product in $(\Lambda_{C_n}^\vee)_{_\QQ}=(\Lambda_{B_n})_{_\QQ}$  
given by the matrix 
%\begin{footnotesize}
	\begin{equation}\label{eqn:scalarBn}
	\left(\vcenter{\xymatrixcolsep{0.0001pc}\xymatrixrowsep{0.0001pc}\xymatrix{
			\scriptstyle 2&\scriptstyle 2&\scriptstyle 2\ar@{.}[rrr]&&&\scriptstyle 2\ar@{.}[rr]&&\scriptstyle 2&\scriptstyle 1\\
			\scriptstyle 2\ar@{.}[dd]&\scriptstyle 4\ar@{.}[dd]&\scriptstyle 4\ar@{.}[rrr]\ar@{.}[dd]&&&\scriptstyle 4\ar@{.}[dd]\ar@{.}[rr]&&\scriptstyle 4\ar@{.}[dd]&\scriptstyle 2\ar@{.}[dd]\\
			&&&&&&&&\\
			\scriptstyle 2\ar@{.}[dd]&\scriptstyle 4\ar@{.}[dd]&\scriptstyle 6\ar@{.}[rrr]\ar@{.}[dd]&&&\scriptstyle 2i\ar@{.}[dd]\ar@{.}[rr]&&\scriptstyle 2i\ar@{.}[dd]&\scriptstyle i\ar@{.}[dd]\\
			&&&&&&&\\
			\scriptstyle 1&\scriptstyle 2&\scriptstyle 3\ar@{.}[rrr]&&&\scriptstyle i\ar@{.}[rr]&&\scriptstyle n-1&\scriptstyle n/2
	}}\right)
	\end{equation}
%\end{footnotesize}

\subsubsection{Case $\# L^c=1$}\ %

\noindent (i) If $L^c=\{1\}$ then  $J_L=\{n\}$, and 
 we deduce from Matrix \eqref{eqn:scalarBn} that
 $n_L=\omega_1^\vee$. Hence,  
$X_{C_n,\{2,\dots,n\}}$ is Gorenstein Fano. Since  
$\#W_L\cdot\omega_n^\vee>n$, it follows that  $X_{C_n,\{2,\dots,n\}}$ is not
smooth. \vspace*{3pt}

\noindent (ii) If $L^c=\{n\}$, then $J_L=\{1\}$ and
$n_L=\omega^\vee_n$. Since  $\Scalar{\omega^\vee_n}{\omega^\vee_i}=i$
if $i\neq n$ and $\Scalar{\omega^\vee_n}{\omega^\vee_n}=n/2$, we
deduce that $X_{C_n,\{1,\dots,n-1\}}$ is Gorenstein Fano if  $n$ is
even, and $\QQ$-Gorenstein Fano of index $2$ if $n$ is odd.

It is easy to see that  $W_L\cdot J_L=W_L\cdot
\omega_1^\vee=\{\omega_1^\vee,\omega_2^\vee-\omega_1^\vee,\ldots,2\omega_{n}^\vee-\omega_{n-1}^\vee\}$
(see  \cite[Remark 5]{kn:voskly}); it follows that  the fan
$\Sigma_{C_n,\{1,\dots,n-1\}}$ is 
simplicial but $X_{C_n,\{1,\dots,n-1\}}$ is not smooth. \vspace*{3pt}

\noindent (iii)  If $L^c=\{2\}$ then  $J_L=\{1,n\}$. Since $
\Scalar{\omega^\vee_1}{\omega^\vee_2}=\Scalar{\omega^\vee_n}{\omega^\vee_2}=2$,
it follows that $n_L=\frac{1}{2}\omega^\vee_2$. By inspection on 
Matrix \eqref{eqn:scalarBn} we deduce that $X_{C_n,\{2\}^c}$ is
Gorenstein Fano.  Since $W_L\cdot J_L>n$, it follows that 
$X_{C_n,\{2\}^c}$ is not smooth.\vspace*{3pt}

\noindent (iv) If $L^c=\{i\}$ with $i\neq 1,2,n$, then
$J_L=\{1,n\}$. Since
$\Scalar{\omega_i^\vee}{\omega^\vee_1-\omega_n^\vee}=2-i$, it
follows that $n_L=0$ and  the variety  $X_{C_n,\{i\}^c}$ is not $\QQ$-Gorenstein Fano.

\subsubsection{Case $M=m+1$}\ %

In this case $J_L=\{1,n\}$.  Assume that
$X_{C_n\{i-1,i\}^c}$, $2\leq i\leq n$ is  $\QQ$-Gorenstein Fano.\vspace*{3pt}

\noindent (i) If  $L ^c=\{n-1,n\}$, then
$n_L=a\omega_{n-1}^\vee+b\omega_n^\vee$, with  $a,b> 0$.  Since  
$n\geq 3$, it follows that
$\Scalar{n_L}{\omega_1^\vee-\omega_n^\vee}= 2a -a(n-1)+b-
\frac{n}{2}b= a(-n+3) -\frac{n-2}{2}b< 0$ and we obtain a
contradiction. \vspace*{3pt}

\noindent (ii) If $L^c=\{i-1,i\}$, $i\neq n$, then 
  $n_L=a\omega_{i-1}^\vee+b\omega_{i}^\vee$,  with  $a,b> 0$.
Since $\Scalar{n_L}{\omega_1^\vee-\omega_n^\vee}=  2a- a(i-1) + 2b -bi=
a(3-i)+ b(2-i)\neq 0$ for all $2\leq i<n$,  we obtain a
contradiction.

\subsubsection{Case $m<M-1$}\ %

In this case $J_L=\{1, m+1,\dots , M-1,n\}$ and,   since $m+1<n$, it
follows that 
$\Scalar{\omega_i^\vee}{\omega_{m+1}^\vee-\omega_1^\vee}\geq 0$ for
all $i\in I$. We deduce from Lemma \ref{lemma:DiscardAlmostAll}
that $X_{C_n,L}$ is not $\QQ$-Gorenstein Fano.

\subsection{Explicit calculations for $G$ of type $D_n$, $n\geq 4$}\ %

Recall that  $D_n\cong(D_n)^\vee$, with  Dynkin diagram
$
\dynkin[labels={1,2,\ ,n-2,n-1,n},label
directions={,,,right,,}] D{}
$. We choose the  $W$-invariant scalar product given by the matrix

%\begin{footnotesize}
	\begin{equation}\label{eqn:scalarDn}
	\left(\vcenter{\xymatrixcolsep{0.0001pc}\xymatrixrowsep{0.0001pc}\xymatrix{
			\scriptstyle 4&\scriptstyle 4&\scriptstyle 4\ar@{.}[rrr]&&&\scriptstyle 4&\scriptstyle 2&\scriptstyle 2\\
			 \scriptstyle 4\ar@{.}[dd]&\scriptstyle 8\ar@{.}[dd]&\scriptstyle 8\ar@{.}[rrr]\ar@{.}[dd]&&&\scriptstyle 8\ar@{.}[dd]&\scriptstyle 4\ar@{.}[dd]&
			\scriptstyle 4\ar@{.}[dd]\\ 
			&&&&&&\\
			\scriptstyle 4\ar@{.}[dd]&\scriptstyle 8\ar@{.}[dd]&\scriptstyle 12\ar@{.}[rrr]\ar@{.}[dd]&&&\scriptstyle 4i\ar@{.}[dd]&\scriptstyle 2i\ar@{.}[dd]
			&\scriptstyle 2i\ar@{.}[dd]\\
			&&&&&&\\
			\scriptstyle 2&\scriptstyle 4&\scriptstyle 6\ar@{.}[rrr]&&&\scriptstyle  2(n-2)&\scriptstyle n&\scriptstyle n- 2\\
\scriptstyle 			2&\scriptstyle 4&\scriptstyle 6\ar@{.}[rrr]&&&\scriptstyle 2(n-2)&\scriptstyle n-2&\scriptstyle n
	}}\right)
	\end{equation}
%\end{footnotesize}

\subsubsection{Case $\#L^c=1$}\ %

\noindent (i)  If $ L^c=\{i\}$, with $i\neq 1,2,n-1,n$, then
$J_L=\{1,n-1,n\}$ and $n_L=a\omega_i^\vee$. Since
$\Scalar{\omega^\vee_i}{\omega^\vee_1-\omega^\vee_n}\neq 
0$, we deduce that $n_L=0$ and   therefore the variety  $X_{D_n,\{i\}^c}$ is not $\QQ$-Gorenstein Fano. \vspace*{3pt}

\noindent (ii) If  $L^c=\{1\}$, then $J_L=\{n-1,n\}$. Clearly
$n_L=\omega^\vee_1/2$, and  $\Scalar{n_L}{\omega^\vee_i}\in \ZZ$ for
all $i\in I$; therefore,  $X_{D_n,\{1\}^c}$ is Gorenstein Fano.
Since  $\# W_L\cdot J_L> n$, it follows that $X_{D_n,\{1\}^c}$ is not smooth
Fano.\vspace*{3pt}

\noindent (iii) If $L^c=\{n\}$, then $J_L=\{1,n-1\}$.
If $n=4$, we deduce by symmetry that $X_{D_4,\{1 \}^c}\cong X_{D_4,\{4 \}^c}$
and therefore   $X_{D_4,\{4\}^c}$ is Gorenstein Fano.

If $n>4$, then
$\Scalar{\omega^\vee_n}{\omega^\vee_1-\omega^\vee_{n-1}}= 4-n>0$. It
follows that  $n_L=0$ and  therefore     $X_{D_n,\{n\}^c}$ is not $\QQ$-Gorenstein Fano. \vspace*{3pt}

\noindent (iv) If $L^c=\{n-1\}$ we deduce by symmetry that
the variety $X_{D_4,\{3\}^c}$ is  Gorenstein Fano, and that
$X_{D_n,\{n-1\}^c}$ is not $\QQ$-Gorenstein Fano if $n>4$.  \vspace*{3pt}

\noindent (v) If $L^c=\{2\}$, then $J_L=\{1,n-1,n\}$. By inspection of
Matrix \eqref{eqn:scalarDn}, we deduce that
$n_L=(\omega^\vee_2)/4$ and it follows that  $X_{D_n,\{2\}^c}$ is
Gorenstein Fano. Again, since   $\# W_L\cdot J_L> n$, if follows that $X_{D_n,\{1\}^c}$ is not
smooth.

\subsubsection{Case $L^c=\{1,2\}$}\ %

 Assume that  $X_{D_n\{1,2\}^c}$ is 
$\QQ$-Gorenstein Fano. Then,   
$n_L=a\omega_1^\vee+b\omega_2^\vee$, with $a,b>0$. Since
$J_L=\{1,n-1,n\}$ and  that  
$\Scalar{n_L}{\omega_n^\vee-\omega_{1}^\vee}= -2a<0$,  we obtain a
contradiction.

\subsubsection{Case $\#L^c>1$, $m>1$}\ %

It is easy to see that  in this case $\{1,n-1,n\}\subset J_L$. Since
$\Scalar{\omega_i^\vee}{\omega_{n-1}^\vee-\omega_1^\vee}\geq 0$ for
$i\geq 2$, it follows from Lemma  \ref{lemma:DiscardAlmostAll} that
$X_{D_n,\{L\}^c}$ is not $\QQ$-Gorenstein Fano.

\subsubsection{Remaining cases}\ %

Since  $\#L^c>1$, $m=1$  and $L^c\neq\{1,2\}$, it follows that
$2<M$; hence, $\{1,2\}\subset J_L=\{ 1,\dots ,M-1, n-1, n\}$. Since
$\Scalar{\omega_i^\vee}{\omega_2^\vee-\omega_1^\vee}\geq 0$, for all
$i\in I$,  it
follows from  Lemma  \ref{lemma:DiscardAlmostAll} that
$X_{D_n,\{L\}^c}$ is not $\QQ$-Gorenstein Fano.

\subsection{Explicit calculations for $G$ of types
	$E_6,E_7,E_8$ and $F_4$}\ %

These cases  can be calculated using Lemma
\ref{lemma:DiscardAlmostAll}  and Theorem
\ref{th:FinalCriteria}, by direct  inspection of the matrices
associated to the $W$-invariant scalar product.  We will
treat in detail only  the cases $E_6$ 
and $F_4$.

\subsubsection{$G$ is of type $E_6$}\ %

The associated Dynkin
diagram of $E_6$ is $\dynkin[labels={1,2,...,6}] E6$. We choose the 
$W$-invariant scalar
product  given by the matrix
\begin{footnotesize}
\[
\left(\begin{matrix}
4 &   3&  5 &   6 & 4 & 2\\
3&    6&    6&    9&    6    &3\\
5&   6& 10 &   12 & 8 & 4\\
6 &   9  &  12 &   18&    12&  6\\
4   & 6 & 8&   12& 10&  5\\
2&    3 & 4  &  6&  5&  4
\end{matrix}
\right)\]
\end{footnotesize}

\noindent (i) If $L^c=\{2\}$, then $J_L=\{1,6\}$, $n_L=\frac{1}{3}\omega^\vee_2$, and
the variety $X_{E_6,\{2\}^c}$  is Gorenstein Fano and not smooth --- because $\Sigma_{E_6,\{2\}^c}$
is not simplicial.\vspace*{3pt}

\noindent (ii) If $L^c=\{i\}$, $i\neq 2$,  then $\#J_L\geq 2$, with
$2\in J_L$, and $n_L=a\omega_i^\vee$. If $j\in J_L\setminus \{2\}$,
then $\Scalar{\omega_i^\vee}{\omega_j^\vee}\neq
\Scalar{\omega_i^\vee}{\omega_2^\vee}$; therefore $n_L=0$ and $X_{E_6,L}$ is not $\QQ$-Gorenstein Fano. \vspace*{3pt}

\noindent (iii) If $L^c>1$, then $\{1,2,6\}\subset J_L$.\vspace*{3pt}

\noindent \emph{(iii-a) $1\notin
L^c$:} since
$\Scalar{\omega_i^\vee}{\omega_2^\vee-\omega_1^\vee}>0$ for all $i\in
L^c$, it follows from Lemma \ref{lemma:DiscardAlmostAll} that $X_{E_6,L} $ is not $\QQ$-Gorenstein Fano.\vspace*{3pt}
% by Lemma \ref{lemma:DiscardAlmostAll}. 

\noindent \emph{(iii-b) $\{1,6\}\subset  L^c$:} since
$\{1,2,3,5,6\}\subset J_L$, it
follows from Lemma \ref{lemma:DiscardAlmostAll} (applied to $j=2$ and
$k=3$) that $X_{E_6,L}$ is not $\QQ$-Gorenstein Fano.\vspace*{3pt}

\noindent \emph{(iii-c) $1\in L^c$ and $6\notin L^c$:} it follows from
Lemma \ref{lemma:DiscardAlmostAll} (applied to $j=2$ and 
$k=6$) that $X_{E_6,L}$ is not $\QQ$-Gorenstein Fano.

\subsubsection{$G$ is of type $E_7,E_8$}\ %

The root systems $E_7$ and $E_8$ have associated  associated Dynkin
diagrams\\
$
\dynkin[labels={1,2,...,7}] E7$  and $ 
\dynkin[labels={1,2,...,8}] E8
$
respectively. 
Notice that in these cases (contrary to the case $E_6$)
we have that 
$\Scalar{\omega_2}{\omega_1}\neq\Scalar{\omega_2}{\omega_n}$, where
$n=7$ or $8$.

% and the scalar products are given by the matrices  
% \[\left(\begin{smallmatrix}
% 4&4&6&8&6&4&2\\
% 4&7&8&12&9&6&3\\
% 6&8&12&16&12&8&4\\
% 8&12&16&24&18&12&6\\
% 6&9&12&18&15&10&5\\
% 4&6&8&12&10&8&4\\
% 2&3&4&6&5&4&3
% \end{smallmatrix}\right)\quad \text{and} \quad
% \left(\begin{smallmatrix}
% 4&5&7&10&8&6&4&2\\
% 5&8&10&15&12&9&6&3\\
% 7&10&14&20&16&12&8&4\\
% 10&15&20&30&24&18&12&6\\
% 8&12&16&24&20&15&10&5\\
% 6&9&12&18&15&12&8&4\\
% 4&6&8&12&10&8&6&3\\
% 2&3&4&6&5&4&3&2
% \end{smallmatrix}\right).\]

We first deal with the case $E_7$; we set $n=7$.

\noindent (i) If $L^c=\{2\}$, then $\{1,n\}= J_L$  and
 $n_L=0$; therefore $E_{n,\{2\}^c}$ is not
$\QQ$-Gorenstein Fano. \vspace*{3pt}

\noindent (ii) If $L^c=\{i\}$, with $i\neq 2$  then
$\#J_L>1$ and $2\in  J_L$ and it
follows as in the case of type $E_6$ (iii-a) that $X_{E_n,\{i\}^c}$ is not
$\QQ$-Gorenstein Fano. \vspace*{3pt}

\noindent (iii)  If $\# L^c>1$, then   either $\{1,2\}\subset J_L$ or
$\{2,n\}\subset J_L$. In both cases  we deduce from Lemma
\ref{lemma:DiscardAlmostAll} that
that $E_{n,L}$ is not $\QQ$-Gorenstein Fano.
\medskip

In order to deal with the case $E_8$, just substitute $n=8$ in the
previous discussion.

\subsubsection{$G$ is of type $F_4$}\ %

% The Dynkin diagrams of $F_4$ and $(F_4)^\vee$ are $\dynkin[label,ordering=Bourbaki]F4$
% and $\dynkin[label,reverse arrows]F4$ respectively. The
% $W$-invariant scalar product in $(\Lambda_R^\vee)_\QQ$ is  
% given by the matrix
The Dynkin diagram of $F_4\cong(F_4)^\vee$ is
$\dynkin[label,ordering=Bourbaki]F4$. We choose the  $W$-invariant scalar product in $(\Lambda_{F_4}^\vee)_\QQ$ 
given by the matrix {\footnotesize $ \left(\begin{matrix}
2&3&4&2\\
3&6&8&4\\
4&8&12&6\\
2&4&6&4
\end{matrix}\right)$}.\vspace*{3pt}

\noindent (i) If $L^c=\{i\}$, with $i\neq 1,4$, then $J_L=\{1,4\}$ and
it follows that $n_L=0$.\vspace*{3pt}

\noindent (ii) If $L^c=\{1\}$, then  $J_L=\{4\}$, so
$n_L=(\omega^\vee_1)/2$. If follows that $X_{F_4,\{2,3,4\}}$ is
$\QQ$-Gorenstein Fano of  Gorenstein index $2$.\vspace*{3pt}

\noindent (iii) If $L^c=\{4\}$, then $J_L=\{1\}$, and $n_L=
\frac{1}{2}\omega^\vee_4$. If follows that $X_{F_4,\{1,2,3\}}$  is Gorenstein Fano but not smooth, since $\# W_L\cdot \omega_1^\vee>4$. \vspace*{3pt}

\noindent (iv) If $\# L^c>1$, then   $\{1,4\}\subset J_L$ and we
deduce from Lemma \ref{lemma:DiscardAlmostAll} that $X_{F_4,L}$ is not
$\QQ$-Gorenstein Fano.  

\section{Some couple of dual reflexive
	polytopes}\label{sec:Reflexive}

% It is well-known (\cite[Theorem 8.3.4]{kn:cox}, see also Proposition \ref{prop:fanogeneral}) that if a toric variety $X_\Sigma$  is Gorenstein Fano, then the polytope $\mathcal P_{R,L}=\Conv(\Prim{\Sigma})$ is reflexive and reciprocally, if $\mathcal P$ is a reflexive polytope, then the toric variety associated to the normal fan of $\mathcal P$ is Gorenstein Fano.  

% If follows that the classification of Gorenstein Fano closures of
% generic orbits  gives a list of couples of reflexive polytopes $\big(\mathcal
% P_{R,L},{\mathcal P_{R,L} }^\circ\big)$.  

Once we have classified all Gorenstein Fano toric varieties $X_{R,L}$, we can   apply \cite[Theorem 8.3.4]{kn:cox} (see
Proposition \ref{prop:fanogeneral}) in order to produce   a list of
couples of reflexive polytopes $\big(\mathcal 
P_{_{R,L}},\mathcal P_{_{R,L}} ^\circ\big)$. Moreover, the $W$-invariance
of $\mathcal P_{_{R,L}}^\circ$ allows us to describe these polytopes as 
Weyl polytopes  (see Definition \ref{def:WeylPolytope}).

% \begin{remark}\label{rem:primsigma-various}
% 	Let $\mathcal P=	\Conv\bigl(\Prim{\Sigma}\bigr)$;  $\mathcal P$ is stable by the
% 	$W$-action, and by duality we deduce that its dual $\mathcal P^\circ\$
% 	is also $W$-stable. Moreover, since $W$ acts transitively on
% 	$\Sigma_{R,L}(n)$, $W$ acts transitively on the maximal proper faces
% 	of $\mathcal P$.

% 	In particular, if $\mathcal Q = \Conv\bigl(\Prim{\sigma_{R,L}}\bigr)$ is $(n-1)$-dimensional with
% 	exterior normal $n_L$, then $\{ w\cdot \mathcal
% 	Q\mathrel{:} w\in W\} $ is the set of proper faces of $\mathcal
% 	P$, and 
% 	if follows  that
%         \marginpar{\color{blue} elimine referencia a \ref{defn:latpol}} 
% 	$
% 	\mathcal P^\circ\= \Conv\bigl( W\cdot (-\varphi_L) \bigr)\subset
% 	(\Lambda_P)_\QQ$ --- recall that  $\Dual{\varphi_{L}}{u}=\Scalar{n_L}{u}$ for all $u\in\Lambda_R^\vee$.
	
% \end{remark} 

\begin{proposition}
	Let $R$ be a root system and $X_{R,L}$
	a Gorenstein Fano generic closure. Then   $\mathcal
	P_{R,L}\subset (\Lambda^\vee_R)_\QQ$ is a dual reflexive
        polytope, with $\mathcal 
	P_{R,L}^\circ = \mathcal{WP}({-\varphi_{L}})\subset 
	(\Lambda_R)_\QQ$,
        where $\varphi_{L}\in\Lambda_R$ is indicated in the last
        column of the corresponding row  of Table \ref{TabFinal}.
\end{proposition}
\begin{proof}
  	It remains to prove the last assertion.
%It is well known that  if $X_{R,L}$ is
%	Gorenstein Fano, then ${\mathcal P_{R,L} =\Conv\bigl(\Prim{\Sigma}\bigr)}$ is% a reflexive
%	polytope.
  Since $\Sigma_{R,L}$ is stable by the
        $W$-action, it follows that $\mathcal P_{R,L}$ and therefore
        $\mathcal P^\circ_{R,L}$ 
       are $W$-stable. Moreover, since $W$ acts transitively on
	$\Sigma_{R,L}(n)$, $W$ acts transitively on the maximal proper faces
	of $\mathcal P_{R,L}$.
In particular,  the set of proper faces of $\mathcal
	P_{R,L}$ is $\{ w\cdot \mathcal
	\Conv\bigl(\Prim{\sigma_{R,L}}\bigr)\mathrel{:} w\in W\} $. Recall that
        $\Dual{\varphi_{L}}{u}=\Scalar{n_L}{u}$ for all
        $u\in\Lambda_R^\vee$; since $\Conv\bigl(\Prim{\sigma_{R,L}}\bigr)$ has $n_L$ as
        exterior normal, it 
	 follows  that
       	\[
	\mathcal P_{R,L}^\circ= \Conv\bigl( W\cdot (-\varphi_L) \bigr)= \mathcal{WP}(-\varphi_L)\subset
	(\Lambda_P)_\QQ.
      \]
\end{proof}   
% If follows from Remark \ref{rem:primsigma-various} that 
% 	$\mathcal P_{R,L}^\vee= \Conv\bigl( W\cdot (-\varphi_L) \bigr)=
% 	\mathcal{WP}(-\varphi_L)$. 

\begin{remark}
	The couple $(L,J_L)$ determines completely the polytope $\mathcal P_{R,L}$.  Indeed, $\mathcal P_{R,L}=\Conv\bigl\langle W\cdot\{\omega^\vee_i\mathrel{ : } i\in J_L\}\bigr\rangle$.  Notice that if $J_L$ contains a single fundamental co-weight $\omega^\vee$, then  $\mathcal P_{R,L}$ is simply the Weyl polytope $\mathcal{WP}({\omega^\vee})$. 	
\end{remark}

We finish this section with an application to the study of root
polytopes --- recall that if $R$ is an (irreducible) root system, then
the \emph{root
	polytope} associated to $R$ is the convex hull
$\Conv(R)\subset (\Lambda_{R})_\QQ$,  see for example \cite{ kn:CelMari1,kn:CelMari2,kn:CelMari3}, where these polytopes are intensively studied.

\begin{proposition}\label{prop:reflpolroot}
	Let $R$ be an irreducible root system. Then the associated root
	polytope $\Conv(R)$ (considered in the root lattice) is reflexive if and only if $R$
	is of type $A_n$,  $C_n$, $D_n$, $E_6$ or $G_2$.

\end{proposition}
\begin{proof}
	Let $\gamma$ be the longest
	root of $R$; then since $\gamma$ is a dominant weight and $-\gamma$ is
	also a root, if follows that $\Conv(R)=\mathcal{WP}({\gamma})=\mathcal{WP}({-\gamma})$. 
	
	Taking into account the description of $\gamma$ as a dominant weight
	(see for example \cite{kn:bourbaki}), we deduce by inspection
        of Table \ref{TabFinal}  that there exists $L$ such that $X_{R,L}$  is
	Gorenstein Fano, with $\varphi_L=\gamma$, if and only if  $R$ is of type $A_n$,  $C_n$, $D_n$,
	$E_6$ or $G_2$.
If this is the case, then	$\Conv(R)=\mathcal{WP}({-\gamma})$ is a reflexive polytope. 
	
	Reciprocally, if $\Conv(R)$ is a reflexive polytope, then the toric
	variety $X$ associated to the normal fan of $\Conv(R)$ is a
	Gorenstein Fano variety. Since  $\Sigma_X=\Sigma_{R,L}$, where
	$\gamma=\sum_{i\in L} a_i\omega_i$, $a_i>0$, the result follows.
\end{proof}

\begin{example} We conclude with three explicit examples of reflexive polytopes
(associated to generic closures) and
their duals.

\noindent \emph{Type $A_2$.} The polytopes
$\mathcal P_{A_2,\{1\}}\cong \mathcal P_{A_2,\{2\}}$  and $\mathcal P_{A_2,\emptyset}$ 
are reflexive. In the figure \ref{picture:Duals}, the
polytopes  $\mathcal 
P_{A_2,L}$, with vertices in 
the weight lattice, are the interior polytopes whereas their duals are
the  exterior ones (with vertex  in the root lattice). 

\begin{figure}[!h]
	%\centering
	\subfloat[$A_2, L=\{1\}$\label{picture:A2_1Dual}]{\resizebox{.29\textwidth}{!}{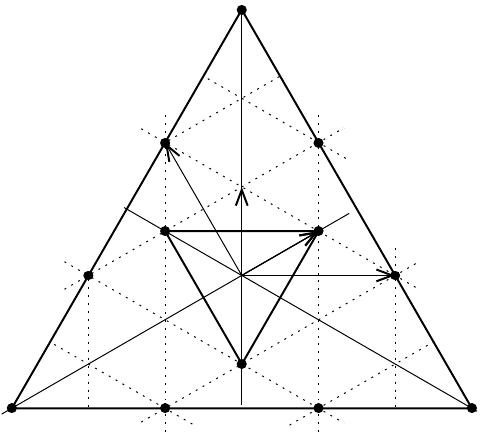}}
	\hspace*{.05\textwidth}
	\subfloat[$A_2,L=\emptyset$\label{picture:A2_12Dual}]{\resizebox{.29\textwidth}{!}{
			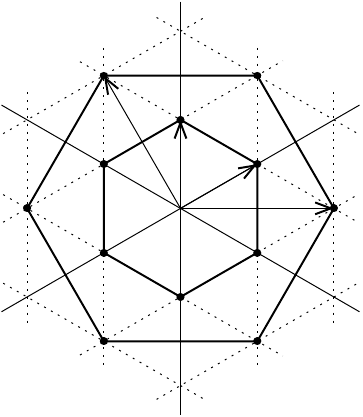}}
	\hspace*{.05\textwidth}	
	\subfloat[$ B_2,L=\{2\}$\label{picture:B2_1Dual}]{\resizebox{.29\textwidth}{!}{
			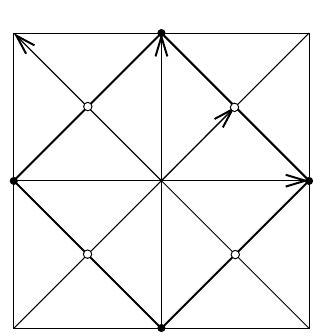}}
	\caption{Three dual reflexive polygons
		associated to root systems of rank $2$. }
	\label{picture:Duals}
	
      \end{figure}
     
\noindent \emph{Type $B_2$}. The polytopes  $\mathcal P_{B_2,\{2\}}$
and $\mathcal P_{B_2,\{1\}}$ are reflexive. In the figure \ref{picture:Duals},  $\mathcal P_{B_2,\{2\}}$ is
the exterior polytope, with vertex  in the weight lattice of $C_2$ ---
it contains 9 lattice points;   $\mathcal P_{B_2,\{2\}}^\vee$  is the
same polytope but in the root lattice of $B_2$, so it contains 5
lattice points.

 \end{example}

% \section*{Declaration}
% Data sharing not applicable to this article as no datasets were generated or analysed during the current study.

\section*{Acknowledgments}

\begin{small} The authors thank the Instituto
Franco-Uruguayo de Matem\'a\-ti\-ca (Uruguay), MathAmSud Project
RepHomol, CSIC (Udelar, Uruguay) and the Institut Montpelli\'erain
Alexander Grothendieck (France) for partial financial support. We
warmly thank C\'edric Bonnaf\'e for his short, but effective introduction
to GAP3. 
\end{small}

 \bibliographystyle{unsrt}
\begin{small}
  \bibliography{BiblioToric}
  \end{small}
\end{document}